\documentclass[11pt]{amsart}
\usepackage{amssymb, amsmath}
\usepackage[table,xcdraw]{xcolor}
\usepackage{ytableau}
 \usepackage{mathtools}
 \usepackage{graphicx}
 \usepackage{pgfplots}
 \usepackage{pifont}
\usepackage{xcolor}
\usepackage{soul}
\usepackage{amsthm}
\usepackage{float}
\usepackage{verbatim} 
\usepackage{multirow}
\usepackage{amssymb}
\usepackage{mathtools}
\usepackage{caption}
\usepackage{subcaption}
\pgfplotsset{compat=1.16}

\definecolor{1}{rgb}{1,0.2,0.3}
\definecolor{2}{rgb}{0.1,0.3,0.5}
\definecolor{3}{rgb}{1,1,0}
\definecolor{4}{rgb}{255,255,255}

\theoremstyle{definition}
\newtheorem{theorem}{Theorem}

\newtheorem{lemma}{Lemma}[section]
\newtheorem{conjecture}[theorem]{Conjecture}

\theoremstyle{remark}

\makeatletter
\@namedef{subjclassname@2020}{\textup{2020} Mathematics Subject Classification}
\makeatother

\begin{document}
\title{Solution Numbers for Eight Blocks to Madness Puzzle}

\author{Inga Johnson}
\address[Inga Johnson]{Department of Mathematics, Willamette University}
\email{ijohnson@willamette.edu}

\author{\'Erika Rold\'an}
\address[\'Erika Rold\'an]{Max Planck Institute for the Mathematics in the Sciences and ScaDS.AI Leipzig University}
\email{roldan@mis.mpg.de}

\thanks{ }

\begin{abstract}
The 30 MacMahon colored cubes have each face painted with one of six colors and every color appears on at least one face. One puzzle involving these cubes is to create a $2\times2\times2$ model with eight distinct MacMahon cubes to recreate a larger version with the external coloring of a specified target cube, also a MacMahon cube, and touching interior faces are the same color. J.H. Conway is credited with arranging the cubes in a $6\times6$ tableau that gives a solution to this puzzle. In fact, the particular set of eight cubes that solves this puzzle can be arranged in exactly \textit{two} distinct ways to solve the puzzle. We study a less restrictive puzzle without requiring interior face matching. We describe solutions to the $2\times2\times2$ puzzle and the number of distinct solutions attainable for a collection of eight cubes. Additionally, given a collection of eight MacMahon cubes, we study the number of target cubes that can be built in a $2\times2\times2$ model.  We calculate the distribution of the number of cubes that can be built over all collections of eight cubes (the maximum number is five) and provide a complete characterization of the collections that can build five distinct cubes.  Furthermore, we identify nine new sets of twelve cubes, called Minimum Universal sets, from which all 30 cubes can be built.

\end{abstract}

\maketitle
\section{Introduction}\label{intro}

\begin{figure}[ht]
    \centering
    \includegraphics[scale=0.50]{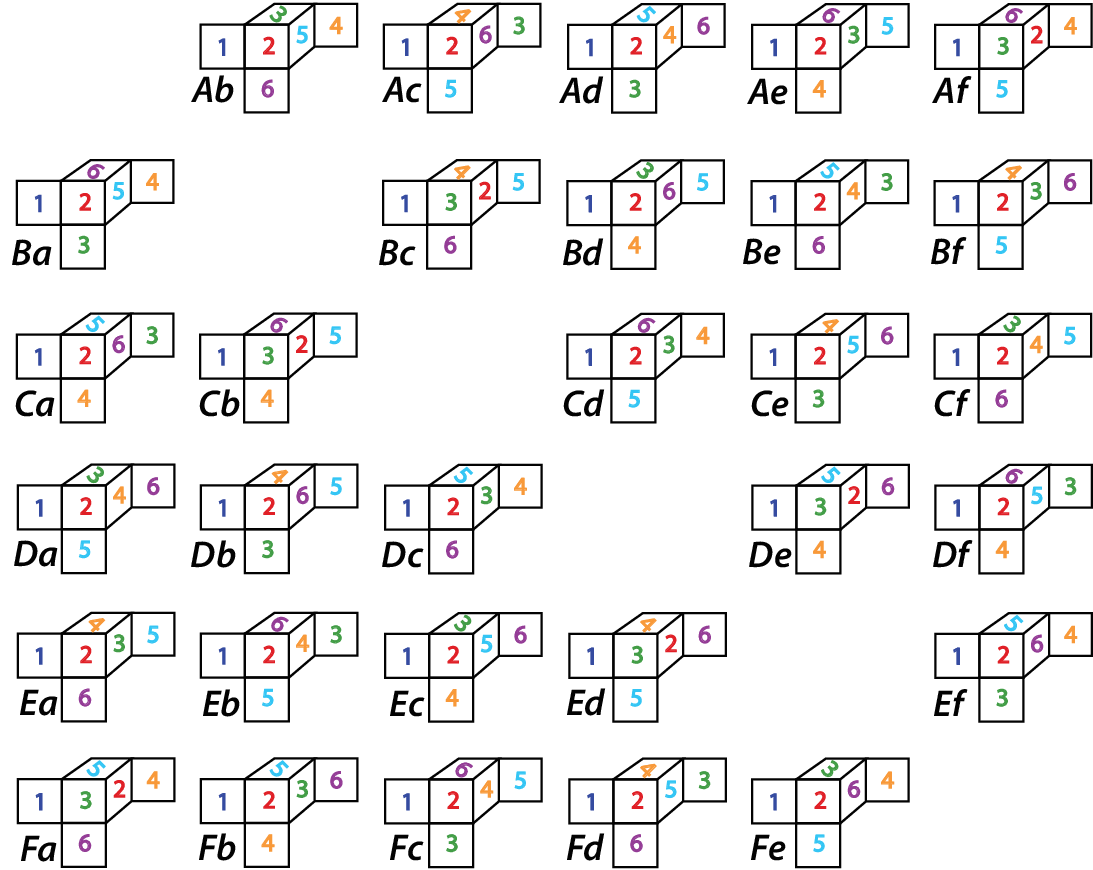}
    \caption{The MacMahon colored cubes with names and tableau arrangement due to J.H.Conway.}
    \label{C6colors}
\end{figure}

The MacMahon colored cubes, also known as 6-colored cubes, are the set of all cubes with each face painted one of six possible colors and every color appearing on at least one face. Rather than working with colors such as red, green, and blue, we designate the six distinct colors with the numbers 1 through 6. Since there are 6!=720 placements of the numbers 1 through 6 on the six faces of a cube and 24 rotational symmetries of a cube, each distinct cube, up to rotational equivalence, appears 24 times in the list of 720 placements. Therefore, there are exactly 30 MacMahon colored cubes as shown in Figure~\ref{C6colors}. 

In his 1921 book \textit{New Mathematical Pastimes}~\cite{MacMahonIntroCubes}, MacMahon introduces several puzzles using the set of 6-colored cubes. The $2\times 2\times 2$ \textbf{interior matching target puzzle} requires stacking eight distinct 6-colored cubes to create a $2\times 2\times 2$ model that is colored as a larger version of a specified \textbf{target cube}, which is also a MacMahon cube, ensuring that touching interior faces in the model have the same color.  MacMahon knew that the interior matching puzzle always has a solution and J.H. Conway is credited with creating the $6 \times 6$ tableau in Figure~\ref{C6colors} that can be used to easily determine the collection of eight cubes that solves this puzzle for a given target.  Interestingly, once the right eight cubes are selected there are exactly {\it two} different ways to arrange them to solve the interior matching puzzle.  We study a less restrictive puzzle, the $2\times 2\times 2$ \textbf{target puzzle}, without the requirement of interior face matching. Our results determine whether a collection of eight distinct 6-colored cubes can be used to solve the puzzle for a given target cube and also determine the number of distinct ways that the eight cubes can be arranged to build the $2\times 2\times 2$ model.  We call this number of ways the \textbf{solution number} of the collection of eight cubes for the target cube.  This $2\times 2\times 2$ \textbf{target puzzle}  was popularized in the early 1970s under the name ``Eight Blocks to Madness" for one collection of eight 6-colored cubes.  Our analysis studies all possible collections of eight colored cubes and applies to all possible targets.  



\subsection{Main Results}\label{sec:Results}

Our results are obtained using computer calculations, modeling and mathematical proofs. All cubes mentioned below are MacMahon colored cubes and we use the term \textbf{collection} to mean a set of eight distinct cubes. There are a total of $\binom{30}{8}=5852925$ collections of colored cubes.

\begin{theorem}\label{solutions}
    Given a target cube,  there are 133680 collections of eight cubes that can be used to build a $2 \times 2 \times 2$ model with the same external coloring as the target. For each of these collections we determine its solution number. In particular, 81  collections achieve the maximum solution number of 16, and 93000 collections have a minimum nonzero solution number of 2. We have also constructed a graph that directly determines the solution number for any collection of eight cubes for a specified target.
\end{theorem}

This means that the probability $p$ of randomly selecting a collection of 8 cubes to build a specific target is around $p \sim 0.0228$. 

\begin{table}[ht]
    \centering
    \begin{tabular}{|c|c|c|c|c|c|c|c|} \hline
       Solution Number  &2&4&6&8&10&12&16 \\ \hline
       Num. of Collections & 93000 & 19860 & 15987 & 2664 & 792 & 1296 & 81 \\
       with that & &&&&&& \\
       Solution Number  & &&&&&& \\ \hline
    \end{tabular}
    \caption{The distribution of the 133680 collections of eight cubes that can be used to build a $2 \times 2 \times 2$ model of a specific target.}
    \label{SolnDist}
\end{table}

This means that the probability $p$ of randomly selecting a collection of 8 cubes to build any target is around $p \sim 0.5259$. 

\begin{theorem}
    Given any collection of eight cubes, at most five targets can be built. Moreover, there are 360 collections that attain this maximum and for each collection none of the five targets are contained in the collection of the eight cubes. 
\end{theorem}\label{thrm:AtMost5Built}

\begin{table}[ht]
\begin{tabular}{|l|c|c|c|c|c|c|} \hline
Number Of &&&&&&\\
Collections & 2774940 & 2256390 & 720405 & 91920 & 8910 & 360  \\ \hline
Number of &&&&&&\\
Targets that can  &&&&&&\\
be built & 0 & 1 & 2 & 3 & 4 & 5  \\ \hline
Proportion of &&&&&&\\
all Collections & .4741 & .3855 & .1231 & .0157 & .0015 & .00006  \\ \hline
\end{tabular}\caption{The distribution of the number of target cubes that can be built over all 5852925 collections of eight MacMahon cubes. It is not possible for a collection of eight colored cubes to build more than five target cubes.}
  \label{DistribBuildable}
\end{table}


A set of colored cubes is called \textbf{universal} if it can be used to build all 30 MacMahon cubes. Haraguchi proves that the minimum size for a universal set is 12 MacMahon cubes and exhibits a set of 12 cubes that can be used to build all 30 target cubes~\cite{haraguchi2016generalization}.  We call Haraguchi's set a \textbf{minimum universal set}.

\begin{theorem}
 There are at least 10 distinct minimum universal sets of 12 cubes. These collections of 12 cubes can build all 30 target cubes.
\end{theorem}
 
\begin{conjecture}
    The 10 minimum universal sets are the only collections of 12 cubes that have this property.
\end{conjecture}

The rest of the paper is structured as follows.
In Section~\ref{sec:background} we provide additional background and a survey of results regarding the MacMahon cubes. In Section~\ref{sec:model}, we describe the graph model used to determine the solution number of a collection for a specific target. In Section~\ref{sec:solutionnumbers}, we prove that our model determines the solution numbers as described in Theorem~\ref{solutions}.
In Section~\ref{sec:all8blocks} we transition to study collections of eight cubes without a specified target in mind. We describe 360 collections of eight cubes that can be used to build five targets and our calculations show that five targets is the maximum number that a collection can build. In Section~\ref{sec:minuniversal} we describe nine new Minimum Universal Sets and some interesting properties of the number of targets that can be built from sets of size $k=8, 9, 10, 11, 12$. In Section~\ref{sec:alg} we provide the algorithm used to calculate results such as those shown in Tables~\ref{SolnDist} and~\ref{DistribBuildable}.


\section{Background}\label{sec:background}
As defined above, the MacMahon colored cubes are cubes colored with six colors. In this section we describe results that extend beyond the MacMahon cubes to include sets of cubes colored with $n$ colors, for $n \in \{2, 3, 4, 5, 6\}$. The set of $n$-colored cubes, denoted by $C_n$, consists of cubes colored so that every face of the cube contains one of the $n$ colors (or a number in $\{1,\dots, n\}$) and each of the $n$ colors appears on at least one face.  

In the 1960s, several puzzles based on colored cubes were built and sold commercially; two of the most popular were {\it 8 Blocks to Madness} and {\it Instant Insanity}.  The {\it 8 Blocks to Madness} puzzle was studied and solved by Kahan~\cite{Kahan} for the specific collection of eight 6-colored cubes in the commercial puzzle.  The {\it Instant Insanity} puzzle contains four distinct cubes from the set of 4-colored cubes with the goal of stacking them vertically so that each $4\times 1$ face of the stack contains each of the four colors.  This puzzle was solved and solution numbers were calculated for any collection of $n$-colored cubes used to build an $n \times 1$ stack for $n \in \{4, 5, 6\}$ by \'E. Rold\'an in~\cite{roldan2018}. Our techniques in this article are inspired by the techniques presented in her paper. Colored cubes are well-known among puzzlers and many other puzzle variations abound \cite{demaine2013variations}. Numerous interesting examples can be found on the website of Jürgen K{\"o}ller ~\cite{KollerWebsite}.

 The \textit{Extended Target Puzzle}, posed by MacMahon~\cite{MacMahonIntroCubes}, is played by selecting a multiset of $k$ $n$-colored cubes with the goal of building a target cube of size $n\times n \times n$, where $k\geq n^3$. We will use the term multiset to refer to sets of colored cubes with possible repetition. Berkove et al. have proven that the $n\times n\times n$ puzzle has a solution if it contains a subset with a $2 \times 2\times 2$ solution. They describe the influence of the action of the symmetric group $S_6$ on the set of 6-colored cubes by way of permuting the colors  in~\cite{Berkove2017AutomorphismsOS}.  The latter paper also studies the minimum cardinality of multisets for which a solution is guaranteed~\cite{berkove2008analysis}. They also describe the influence of the action of the symmetric group $S_6$ on the set of 6-colored cubes by way of permuting the colors.  Berkove and collaborators proceed to study the extended target puzzle for  2-, 3-, and 4-colored cubes~\cite{berkove2017color}. From the set of 4-colored cubes, they provide a collection of 10 cubes from which all 68 targets can be built~\cite{berkove2018michael}. More recently, they have studied the extended target puzzle with 2-colored hypercubes~\cite{newBerkove}.

Building on Berkove's work, Haraguchi studied several optimization questions regarding the MacMahon cubes~\cite{haraguchi2016generalization}.  He identified a multiset of 23 cubes from which no solution to the $2\times 2 \times 2$ target puzzle exists and proved this to be the largest size for such a set. These 23 cubes are referred to as a \textbf{maximum infeasible} (multi)set. He also discovered a subset of 12 distinct MacMahon cubes for which the $2\times 2 \times 2$ target puzzle can be solved for all 30 of the 6-colored cubes and proved that 12 is the minimum size of a set with this universal property.

\section{A Model for the Target Puzzle}\label{sec:model}

 Using the notation of  J. H. Conway as shown in Figure~\ref{C6colors}, each MacMahon cube is described by a \textbf{two letter name}. The first letter, one of A, B, C, D, E, or F, indicates the row position of the cube in the tableau and the second letter a, b, c, d, e, or f, indicates the column position. Additionally, we denote each cube in the tableau by a collection of eight corner numbers. A \textbf{corner number} is the three-digit number of least numerical value that cyclically describes the numbers around a cube corner when read in the clockwise direction.  The corner numbers for Fb are shown in Figure~\ref{exC36_1}. Each MacMahon cube and its eight corner numbers can be found in Table~\ref{C36corners} of the Appendix. 

\begin{figure}[h]
    \centering
    \includegraphics[scale=0.8]{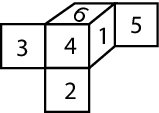}
     \caption{The eight corner numbers of the cube Fb are $\{ 124, 146, 165, 152, 234, 253, 356, 364\}$.}
    \label{exC36_1}
\end{figure} 

We observe that there are exactly 40 possible corner numbers and within every row of Conway's tableau all 40 corner numbers are present exactly once. The same holds for any column of the tableau.  Referring back to Figure~\ref{C6colors}, notice that for any cube, Xy, the cube Yx, contains the same color positions on the vertical sides but the colors on the top and bottom of the cube are swapped.  These pairs, Xy and Yx, are called mirror cubes and they contain no corner numbers in common. For example, if Xy contains 123, then Yx contains 132. In general, for any target cube, Xy, the eight cubes that solve the interior matching puzzle are precisely those that are in the same row or same column as the target's mirror cube, Yx.

We begin our analysis by considering a specific target, for example Ba, and any collection of eight cubes.  Our focus on this one target is sufficient since a bijective recoloring will map solutions for any other target to the corresponding solutions found here for Ba.  

The \textbf{usable corner count} for a colored cube, C, with respect to a given target, T, is the number of corner numbers of C that are also corner numbers of T.  For the MacMahon colored cubes the usable corner count of a cube for any target is either 0, 2, or 8, where 8 only arises when C=T.  When a cube has a usable corner count of 0 for the target cube, we say that it is \textbf{unusable} in the target puzzle, since any collection of eight cubes with one or more unusable cubes will produce no solutions for the target T. 

\begin{figure}[H]
    \centering \includegraphics[scale=0.2]{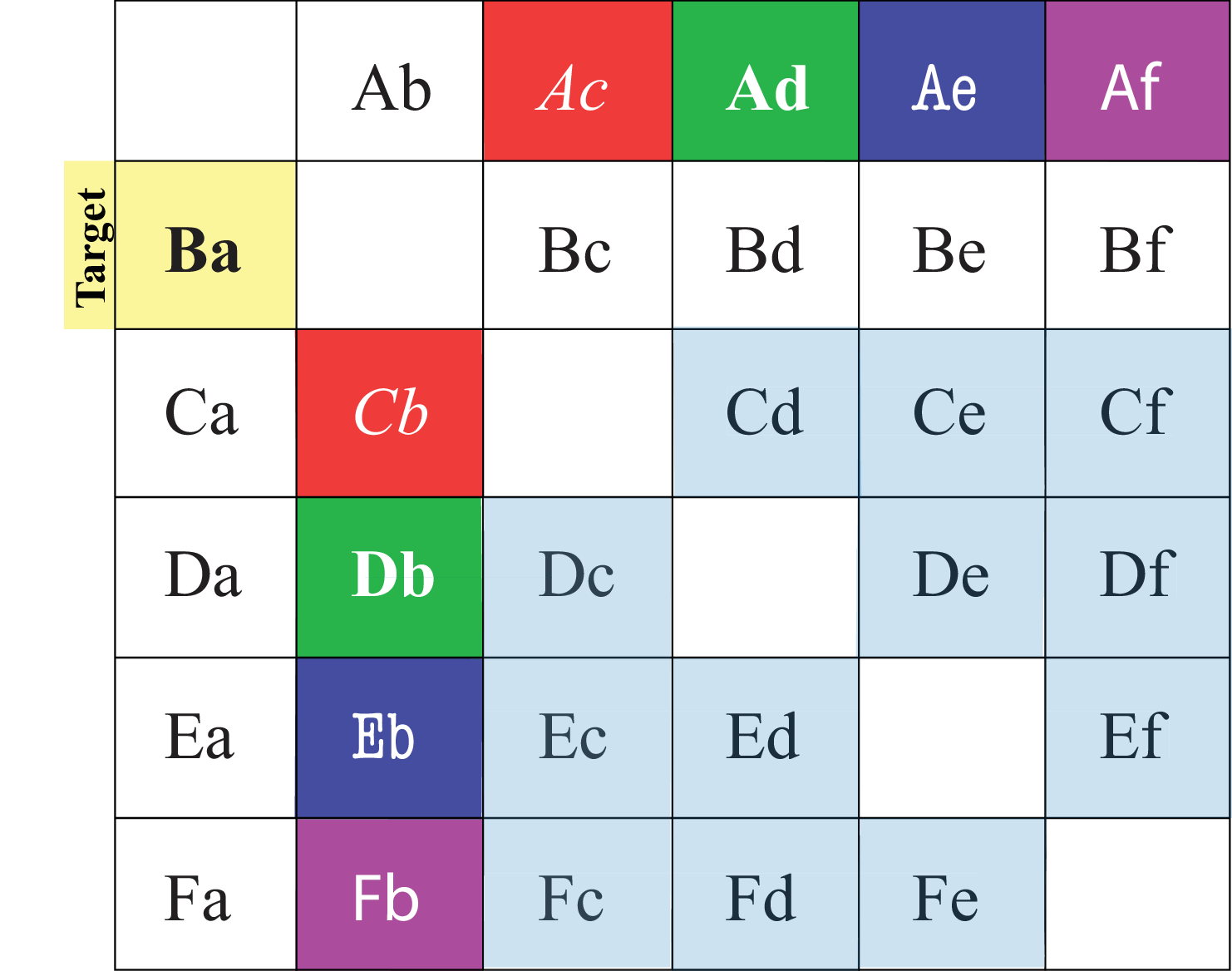}
\caption{The usable cubes for target Ba are highlighted with a shaded background in Conway's tableau. The color/font pairs  $\textit{Ac, Cb}$; $\textbf{Ad, Db}$; $\texttt{ Ae, Eb}$; and {\fontfamily{lmss}\selectfont Af, Fb} indicate  cubes pairs that share the same two corner numbers with the target Ba.}
    \label{IntMatchSoln}
\end{figure}

Nine cubes are unusable for our target cube Ba: the cubes in the same row Bc, Bd, Be, Bf; the cubes in the same column Ca, Da, Ea, Fa; and the mirror cube Ab. The target cube has a usable corner count of 8, and the remaining 20 cubes have exactly two corners in common with the target cube Ba.  Figure~\ref{IntMatchSoln} depicts, with a shaded background, the 21 cubes with at least one corner number in common with the target. We refer to these 21 cubes as \textbf{usable cubes} for the Ba target puzzle.  The cubes Ac, Ad, Ae, Af, Cb, Db, Eb, Fb are the eight cubes that can be used to solve the $2\times 2\times2$ puzzle with interior matching in two distinct ways.   The color and font pairings Ac, Cb; Ad, Db; Ae, Eb; and Af, Fb; indicate that the two cubes share the same two corner numbers in common with the target Ba.  

\begin{figure}
    \centering
    \includegraphics[scale=0.2]{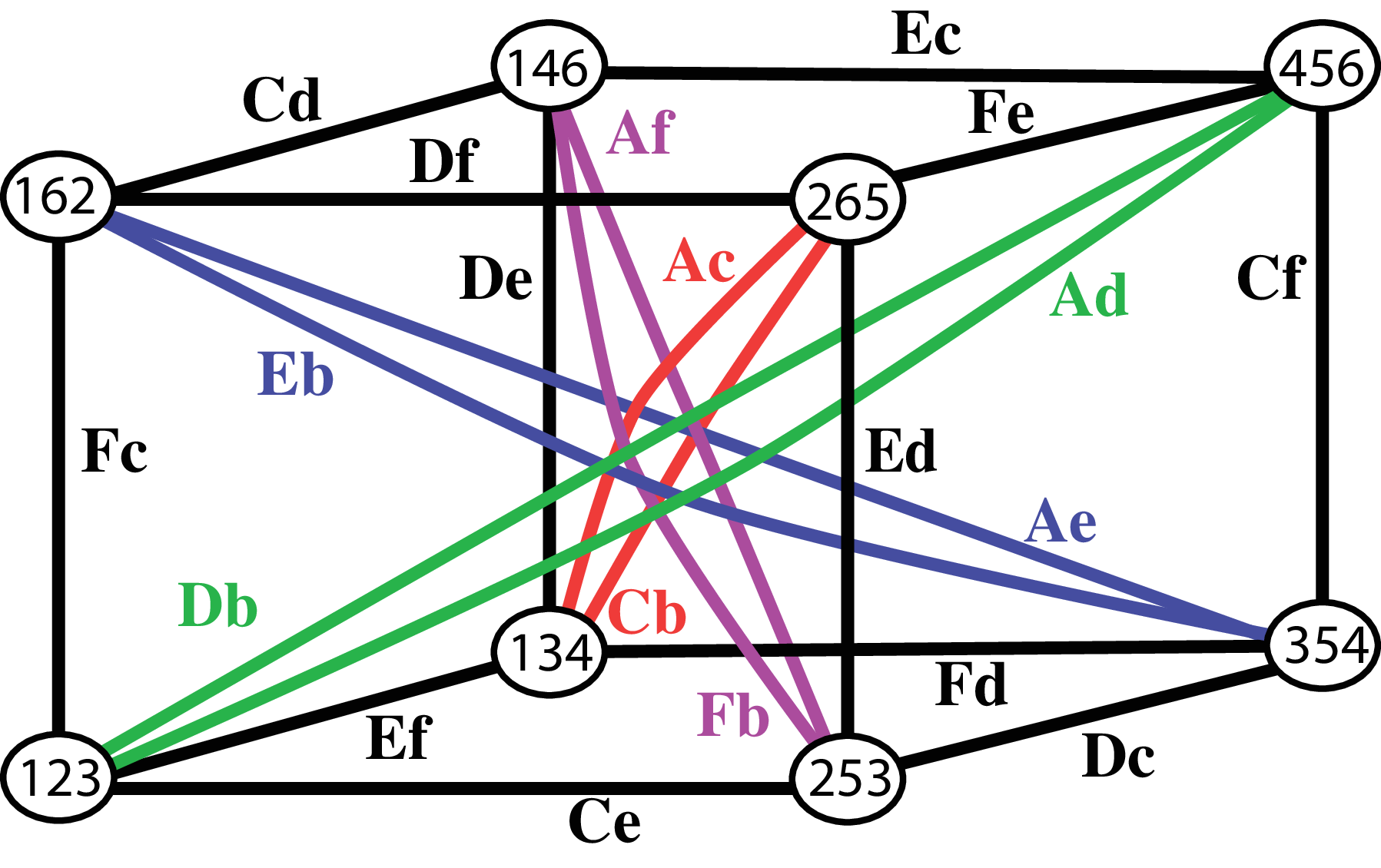}
    \caption{A geometric representation of the graph used to calculate solution numbers.}
    \label{GeometricRepBa}
\end{figure}

Our model depicts the 20 cubes with two corners in common with Ba as edges in a graph, as shown in Figure~\ref{GeometricRepBa}. This graph is a geometric representation of the cube Ba, with its eight corner numbers shown as labeled vertices.  The edges in the graph are labeled with the name of the cube that shares with the target the two corner numbers (vertices) that the edge connects.  We observe that the 12 standard edges in the geometric representation of a cube are labeled by cube names Cd, Ce, Cf, Dc, De, Df, Ec, Ed, Ef, Fc, Fd, and Fe. The remaining eight cubes, which are the same eight cubes that form the solution to the interior matching puzzle, share two diagonally opposite corner numbers across the center of the cube. The only usable cube for the Ba target puzzle that is not depicted in Figure~\ref{GeometricRepBa} is the target Ba itself. The paired cube names and edge colors in Figure~\ref{GeometricRepBa} align with the color pairs shown in Figure~\ref{IntMatchSoln}. 

The graph in Figure~\ref{GeometricRepBa} is crucial for determining the solution number for a collection of eight cubes. Let $M$ is the graph in Figure~\ref{GeometricRepBa} and $W$ be a collection of eight cubes.  Let $G_W$ denote the subgraph of $M$ with edges corresponding to the cubes in $W$ and all eight vertices (corner numbers). If the subgraph $G_W$ contains fewer than seven edges,  the collection $W$ yields no solution for the target Ba, since at least one cube in $W$ is unusable.

\begin{figure}
    \centering
    \includegraphics[scale=0.2]{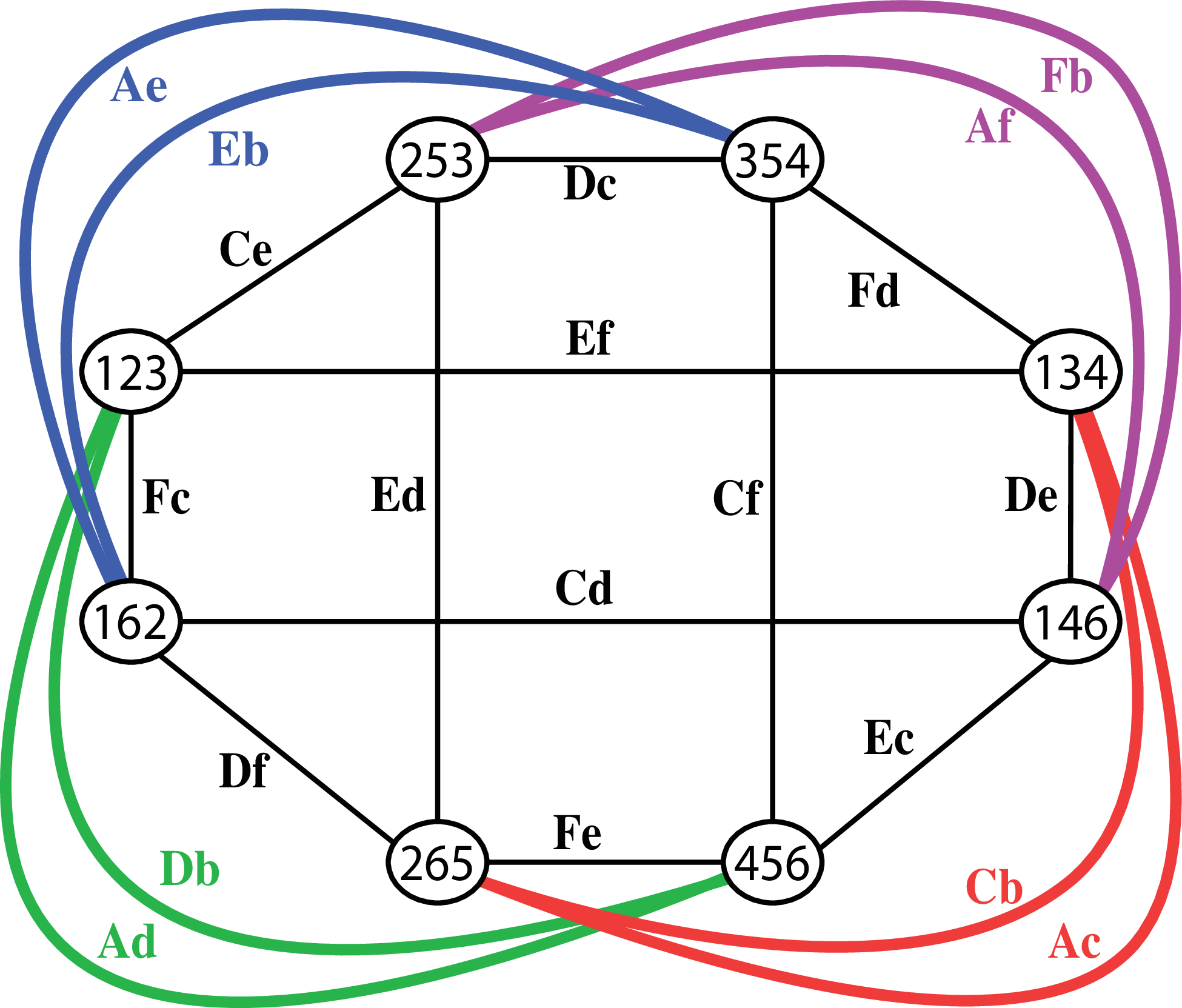}
    \caption{The graph used to determine the solution number for the target Ba and any collection of eight MacMahon cubes.}
    \label{BaCornerAltView}
\end{figure}

If the subgraph $G_W$ contains 8 edges and a vertex of degree 0, then none of the cubes in the collection contain the corner number corresponding to the vertex of degree 0. Consequently, no cube in $W$ can fill that corner of the $2 \times 2 \times 2$ target puzzle, resulting in a solution number of zero.  In other cases, we determine the solution number of the collection $W$ by examining the number of components, cycles, and trees in the subgraph $G_W$.  The cycles in $M$ are easier to see in Figure~\ref{BaCornerAltView}, which is an isomorphic representation of the graph in Figure~\ref{GeometricRepBa}. We continue to refer to Ac, Ad, Ae, Af, Cb, Db, Eb, and Fb as `diagonals' for Ba.  We share a few examples before proving Theorem 1.

Consider the collection $W$ containing the eight cubes that form the solution for the interior matching puzzle: Ac, Ad, Ae, Af, Cb, Db, Eb, and Fb. This collection of eight cubes results in a subgraph $G_W$ containing the four disjoint 2-cycles from the graph $M$.  Each 2-cycle consists of a pair of cubes with the same two corners in common with Ba. Each cube in a pair has two possible placements; for example, Ae can be placed in corner 162 and Eb in corner 354, or vice versa. Thus, each cycle in $G_W$ contributes a factor of 2 to the solution number; hence, the solution number for this collection with four cycles is $2^4 =16$.  

Furthermore, any one of these eight cubes could be replaced by the target cube Ba. Any such collection will also have a solution number of 16 because the corner numbers of the replaced cube can and must be satisfied by the target cube. Thus, we have described a total of 9 collections that all achieve the maximum number of 16 solutions.  

There are an additional 72 collections of eight cubes that attain the maximum solution number of 16 for our target Ba. These additional 72 collections all contain the target Ba, and the remaining seven cubes are selected as: two of the 2-cycle pairs, and a path of three cubes that spans the four vertices not contained in the two 2-cycles. There are 6 ways to select two of the four 2-cycles, and a straightforward count shows that there are 12 distinct paths of length 3 that span the four remaining vertices.  Therefore, we have 6 ways to pick the 2-cycle pairs and 12 ways to select a path of three cubes connecting the remaining four corner numbers, resulting in $72$ collections of eight cubes.  

Each of these 72 collections has a solution number of 16 because the two 2-cycles each contribute a factor of 2 to the solution number, and the target cube can be placed in any one of the four corner numbers along the path of the three cubes forcing the other three cubes along the path to be placed in the only remaining possible spot. The path of length three (with 4 vertices) and the target cube contribute a factor of 4 to the solution number.  Thus, these 72 collections all yield 16 distinct solutions for the Ba target puzzle.  Combined with the 9 collections mentioned above, this accounts for all  81 collections of 8 cubes with a solution number of 16 as reported in Table~\ref{SolnDist}. The algorithm in Section~\ref{sec:alg} was used to calculate the remaining entries of the table, and Theorem~\ref{solutions} can be used to verify the counts.
 
\section{Solution Numbers and a Proof of Theorem 1}\label{sec:solutionnumbers}

Let $W$ be a collection of eight MacMahon cubes. Let $M$ denote the graph in Figure~\ref{GeometricRepBa}, or equivalently Figure~\ref{BaCornerAltView}. Let $G_W$ be the subgraph of $M$ that contains all edges labeled with cube names from $W$ and  all vertices of $M$.  We will prove that the solution number of $W$ for the target Ba is completely determined by $G_W$ and whether or not the target cube, Ba, is in $W$. All arguments are intertwined with the graph $M$, its properties, and its possible subgraphs.

The following lemma describes several cases where the solution number of a collection of eight cubes is zero.

\begin{lemma}\label{noSolution}
   If the graph $G_W$ satisfies any of the following properties, 
   \begin{itemize}
       \item $G_W$ contains fewer than eight edges and the target Ba is not in $W$; 
       \item $G_W$ contains a component that is a tree and the target Ba is not in $W$; or
       \item $G_W$ contains more than one component that is a tree, 
   \end{itemize}  
   then the solution number of $W$ is 0.
\end{lemma}

\begin{proof}  Suppose $G_W$ contains fewer than eight edges and the target Ba is not in $W$. Since $W$ doesn't contain the target cube, it must be the case that $W$ contains a cube with no usable corner numbers for the target cube. Therefore $W$ has a solution number of 0 for the target Ba. 

When the graph $G_W$ contains a component that is a tree with $k$ edges, then the $k$ cubes along the edges are insufficient to satisfy the $k+1$ corner numbers in the tree. For such a collection to have a solution, the target cube must be used to satisfy one of the corner numbers in the tree. If we root the tree at the corner number where the target is placed, then tree edge cubes have forced placements in the corner numbers away from the root.  So, a collection without Ba and a tree in $G_W$ results in no solution, as does a collection with more than one tree in $G_W$.
\end{proof}

Note that if $G_W$ contains eight edges in $M$, then the target cube is not in $W$.   Theorem~\ref{WnoTarget} describes the solution number for collections that do not contain the target cube $Ba$, and Theorem~\ref{WhasTarget} describes the solution number for collections that contain the target. Together, these results prove Theorem 1.

\begin{theorem}\label{WnoTarget}
    Suppose $G_W$ is a subgraph of $M$ containing eight edges.  Let $n$ be the number of connected components in $G_W$. If none of the $n$ connected components are trees, then the solution number of $W$ for the target Ba is $2^n$. Otherwise, the solution number is 0.
\end{theorem}

\begin{theorem}\label{WhasTarget}
Suppose the target cube Ba is in $W$ and $G_W$ is a subgraph of $M$ containing seven edges. If $G_W$ contains exactly one connected component that is a tree with $k$ edges (so all other connected components contain cycles), then the solution number of $W$ for the target is given in the table below.
\begin{center}
   \begin{tabular}{|c|c|c|c|c|} \hline
        $\#$ of connected components in $G_W$ & 1 & 2 & 3 & 4 \\ \hline
        solution number of $W$ & 8 & 2(k+1) & 4(k+1) & 16 \\ \hline
    \end{tabular}
    \end{center}
If $G_W$ has more than one component that is a tree, then by Lemma~\ref{noSolution} there are no solutions. 
\end{theorem}

We prove Theorem~\ref{WnoTarget} through a sequence of lemmas.

\begin{lemma}\label{noTarget1component}
Suppose $G_W$ is a subgraph of $M$ containing eight edges. If $G_W$ is a connected graph, then the solution number of $W$ is 2. \end{lemma}

\begin{proof} Suppose $G_W$ is a connected graph that spans all eight vertices of $M$. Then $G_W$ has a subgraph that is a spanning tree with seven edges and all eight vertices. Therefore the eighth edge, when added to the spanning tree, will form a single cycle in $G_W$. To determine a solution to the puzzle we orient the edges of $G_W$ to indicate the corner number where each cube (a labeled edge) will be placed in the corners of Ba (a vertex of $M$).  One such orientation of the edges is to point each edge of the cycle in the direction of traversing around the cycle. There are two possible ways to traverse the cycle, which leads to two distinct placements of the set of cubes along the cycle.  Since the graph $G_W$ is connected, all remaining edges are connected to the cycle by some path and can be oriented to point away from the cycle.  The cubes not in the cycle have a forced placement in the puzzle solution. The graph contains all eight vertices, so every corner number will be satisfied by the cube on the edge pointing towards it. Thus, the collection $W$ has a solution number of 2. \end{proof}

\begin{lemma}\label{noTarget2}
Suppose $G_W$ is a subgraph of $M$ containing eight edges. If $G_W$ consists of two connected components, neither of which are trees, then the solution number of $W$ is 4. 
\end{lemma}


\begin{proof} Suppose $G_W$ consists of two connected components.  The eight vertices must be split between the two connected components as 1 and 7, 2 and 6 , 3 and 5, or 4 and 4 vertices.  The case with 1 and 7 vertices cannot happen, as a single vertex component that is not a tree would imply the existence of a loop, but the graph $M$ contains no loops. For the three remaining cases, both components contain a spanning tree. Therefore, the respective spanning trees must contain a total of 6 edges that are split between the two components as 1 and 5 edges,  2 and 4 edges, or both spanning trees have 3 edges. Since $G_W$ has eight edges, two of its edges are not in either spanning tree. If both of these two edges are in the same connected component of $G_W$ then the other component is a tree, which contradicts our hypotheses. Thus each connected component contains an edge that is not in the spanning tree, and so each connected component contains one cycle. By the same argument used in Lemma~\ref{noTarget1component}, a component in $G_W$ that contains one cycle contributes a factor of 2 to the solution number. Thus, the solution number of $W$ is $4$. \end{proof}

\begin{lemma}\label{noTarget3or4}
Suppose $G_W$ is a subgraph of $M$ containing eight edges. If $G_W$ consists of three or four connected components, none of which are trees, then the solution number of $W$ is 8 or 16, respectively. 
\end{lemma}

The proof of Lemma~\ref{noTarget3or4} is similar to the proof of Lemma~\ref{noTarget2}.  The solution number increases by a factor of two for each component with a single cycle. It cannot be the case that a component has more than one cycle because some other component would then be a tree.
Lemmas~\ref{noTarget1component}, \ref{noTarget2} and~\ref{noTarget3or4} complete the proof of Theorem~\ref{WhasTarget}.

Next we prove Theorem~\ref{WhasTarget} through a sequence of lemmas about collections that contain the target cube Ba.

\begin{lemma}\label{wTarget1component}
Suppose the target cube Ba is in $W$ and $G_W$ is a subgraph of $M$ containing seven edges. If $G_W$ is a connected graph, then the solution number of $W$ is 8. 
\end{lemma}

\begin{proof} Suppose $G_W$ has seven edges, all eight vertices, and is connected. Then $G_W$ is a tree.  The target cube can satisfy any of the eight corner numbers. We place the target cube at any vertex, then root the tree at that vertex and each of the seven cubes in the collection satisfies a unique corner number by directing the edges away from the root. Thus the solution number for $W$ is 8, corresponding to the number of possible placements of Ba. \end{proof}

\begin{lemma}\label{wTarget2components}
Suppose the target cube Ba is in $W$ and $G_W$ is a subgraph of $M$ containing seven edges. If $G_W$ is disconnected with two components, then one of the components is a tree with $k$ edges,  for $k\in \{0, 1, 2, 3, 4, 5\}$, and the solution number of $W$ is $2(k+1)$. 
\end{lemma} 

\begin{proof}
Suppose $G_W$ contains seven edges and two components.  The components will contain one of the following combinations of edges; 7 and 0 edges, 6 and 1 edge, 5 and 2 edges, 4 and 3 edges. We consider each case:

Case 1 (7 and 0 edges): The 0 edged component is a single vertex, so a tree.  The component with 7 edges must contain the other 7 vertices which means this component contains one cycle. In this case the target cube Ba satisfies the corner number of the single vertex component and the component with 7 edges contributes a factor of 2 to the solution number. 

Case 2 (6 and 1 edges):  The component with one edge is a tree and includes two vertices. The remaining 6 vertices are in the other component with 6 edges; therefore this component contains one cycle.  Thus the component with 6 edges contributes a factor of two to the solution number. Since the target cube is in $W$, it can and must be used to satisfy one of the two corner numbers in the component with one edge.  So the target cube also contributes a factor of two to the solution number.  Thus, the solution number is 4. 

Case 3 (5 and 2 edges): Either the 2-edged component is a cycle containing 2 vertices or a tree with 3 vertices. Assuming the case of the cycle with two vertices, the remaining 6 vertices must be in the component with 5 edges. Thus that connected component is a tree. Similar to the argument in the case above, to create a solution the target cube can and must be placed in one of the 6 vertices in the tree, resulting in a factor of 6 in the solution number. The cycle of length 2 gives another factor of two, so the solution number is $12$ when $k=5$. 

Else, the component with 2 edges is a tree with 3 vertices.  In this case, the the other 5 vertices are in the component with 5 edges and thus contain one cycle. The cycle contributes a factor of 2 to the solution number and the target cube can be used in any of the three corner numbers in the tree component, resulting in a factor of 3 in the solution number. So the solution number of the $2\times2\times2$ puzzle for the collection $W$ is 6  when $k=2$. 

Case 4 (4 and 3 edges): The 3-edged component is either a tree or contains a cycle.  If the three edge component contains a cycle, it must contain exactly one cycle of length 2, since the graph is bipartite and has no odd cycles. Thus the 3 edged component contains three vertices.  Then the remaining 5 vertices must all be in the connected component containing 4 edges, and thus that component is a tree. So in this case, the solution number of $W$ is 10 when $k=4$. 
Else, the 3-edged component is a tree and contains 4 vertices.  This implies the remaining 4 vertices are in the connected component with 4 edges which forces this component to contain one cycle. So in this case, the solution number is $8$ when $k=3$.

Notice that $k$ cannot be greater than or equal to 6 because a tree of 6 edges would include 7 vertices, leaving one vertex for the second component. However, the graph $M$ contains no loops, so the one vertex must have degree 0.  This explains why a solution number of 14 is not possible. If it were possible, it would come from a cycle that contributes a factor of 2  and a tree with 7 vertices. Since the cycle would require 2 vertices and the tree would need 7 vertices, this case is not possible in a disconnected graph with a total of 8 vertices. \end{proof}

\begin{lemma}\label{wTarget3components}
Suppose the target cube Ba is in $W$ and $G_W$ is a subgraph of $M$ containing seven edges and all eight vertices. If $G_W$ is disconnected with three components and  exactly one of the three components is a tree with $k$ edges,  where $k\in \{1, 2, 3\}$, then the solution number of $W$ is $4(k+1)$.  Otherwise, the solution number is 0. \end{lemma}

\begin{proof} 
Suppose $G_W$ is a subgraph of $M$ containing seven edges and all eight vertices and consists of three connected components. For a component to possibly contribute to a nonzero solution number, it must contain at least one cube.  Thus, each component must have at least one edge. Therefore, the seven edges could be distributed into three components in the following four ways:
    1, 1, 5; \hspace{0.1in} 1, 2, 4; \hspace{0.1in}1, 3, 3; \hspace{0.1in} or  2, 2, 3.

We will consider each case separately.

Case 1 (components contain 1, 1, and 5 edges): The two components with one edge each imply that graph $G_W$ contains more than one tree. This results in a solution number of 0 by Lemma~\ref{noSolution}.

Case 2 (components contain 1, 2, and 4 edges): Since there is only one tree,  the component with one edge must be a tree with 2 vertices. For a nonzero solution,  the components with 2 and 4 edges cannot be trees.  Therefore the component with 2 edges is a 2-cycle with 2 vertices, and the component with 4 edges must have the remaining 4 vertices and thus a cycle. Each component with a cycle contributes a factor of 2 to the solution number and the tree with 2 vertices also contributes a factor of 2.  Thus, in this case with $k=1$, the solution number is $8$.

Case 3 (components contain 1, 3, and 3 edges):  As above, the component with one edge must be the tree with 2 vertices. The components with 3 edges cannot be trees and must have at least 3 vertices. In fact, the only possible vertex count for both of these two components is exactly 3 vertices. Thus, each of these two components must have a cycle. As above, $k=1$ and the solution number is $8$.

Case 4 (components contain 2, 2, and 3 edges): The tree can either contain 2 or 3 edges.  If one of the components with 2 edges is the tree, then then that component contains 3 vertices. The component with 3 edges must contain at least 3 vertices, leaving at most 2 vertices for the other component with 2 edges.   So we have a tree with 3 vertices, and the other two components must contain a cycle. Therefore, with $k=2$, the solution number for $W$ is $12$.  
Alternatively, if the component with 3 edges is the tree, then $G_W$ contains a tree with 3 edges and 4 vertices.  Consequently, the other two components, each with 2 edges, must both have exactly 2 vertices.  As a result, these two components are both cycles, and the solution number is 16. \end{proof}

The last case to consider is when $G_W$ contains 4 connected components. 

\begin{lemma}\label{wTfourcomponents}
Suppose the target cube Ba is in $W$ and $G_W$ is a subgraph of $M$ containing seven edges. If $G_W$ is disconnected with four components, then the solution number of $W$ is 16 if there is only one tree, otherwise it is 0. 
\end{lemma}

\begin{proof}
From Lemma~\ref{noSolution}, a subgraph with two or more trees results in a solution number of 0.  So we consider the case where there is only one tree. In this case, the edges must be distributed to the four components as 1, 2, 2, and 2 edges. The component with one edge must be the tree, while each component with 2 edges is a 2-cycle. The target cube and the tree contribute a factor of 2 to the solution number, as does each of the three 2-cycle components. Thus, the solution number is 16. \end{proof}

Lemmas~\ref{wTarget1component}, \ref{wTarget2components}, \ref{wTarget3components}, and \ref{wTfourcomponents} complete the proof of Theorem~\ref{WhasTarget}.

\section{All \textit{Eight Blocks to Madness} puzzles and a Proof of Theorem 2}\label{sec:all8blocks}

In this section, we pivot from calculating solution numbers to considering all $5852925$ collections (30 choose 8) and investigate the number of target cubes that each collection can build.  Our calculations show that $2774940$ collections can build no target cubes, and the remaining collections can build 2, 3, 4 or 5 cubes, as shown in Figure~\ref{thrm:AtMost5Built}. There are no collections that build more than 5 target cubes.
Our calculations show that an elusive 360 collections can be used to build 5 target cubes.  What do these collections look like and which sets of five targets can they build?  We answer this question by constructing a set of five targets and the collection that can build them according to the rules in Theorem~\ref{5targets}.

\begin{theorem}\label{5targets} A collection of eight cubes constructed as follows can be used to build five target cubes.  The five targets are specified below.

 \textbf{The five targets:} Start by selecting exactly three columns of Conway's tableau and refer to the columns as $P:=\{c_1, c_2, c_3\}\subset\{a, b, c, d, e, f\}$. Select four rows of the tableau so that exactly two of them correspond to selected column letters in $P$ and the other two rows are not in $P$. Now consider the complement of the union of the 4 rows and 3 columns selected above; this set will contain five cubes and the remaining blank diagonal entry of the tableau. These five cubes are the targets that can be built.  

\textbf{The collection of 8 cubes:} There are 10 cubes in the intersection of the selected 4 rows and 3 columns.  Of these 10 cubes, exactly two will be mirror cubes of two of the five target cubes. The other eight cubes comprise our collection. 

Analogously, we can select three rows followed by four columns so that exactly two columns correspond to selected rows and two columns do not. 
\end{theorem}

\begin{proof}  There are 20 ways to select the three columns and $9$ ways to select the four rows; thus, following the selection rules results in 180 collections. Analogous rules to select three rows followed by four columns results in another 180 collections.  This gives a total of 360 collections, which aligns with the number of collections found using our algorithm. Our calculations show that these are exactly the 360 collections that can build the five specified target cubes.
\end{proof}

\begin{figure}[h]
   \centering  \includegraphics[scale=0.2]{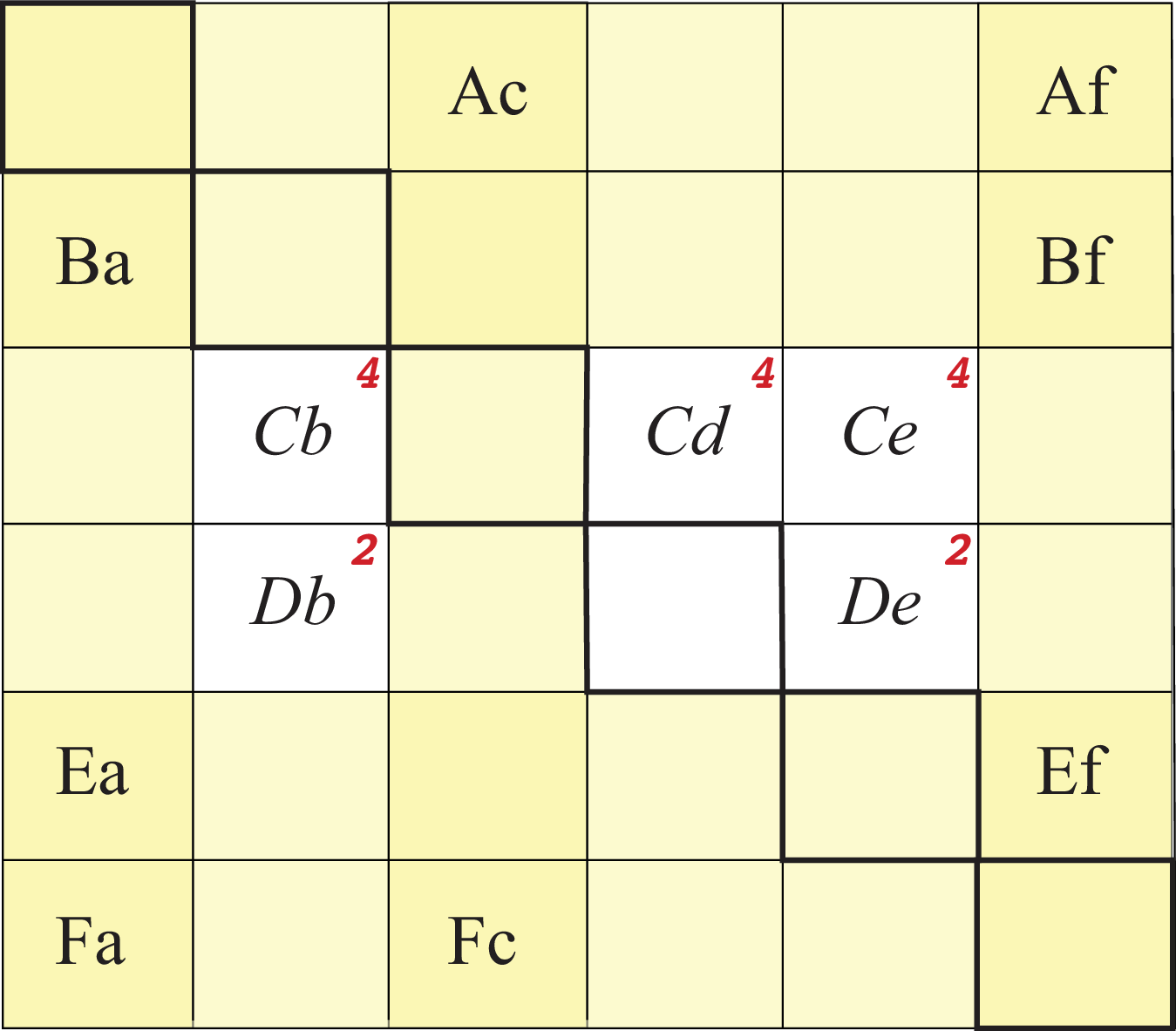} 
   \caption{The collection {\bf Ac, Af, Ba, Bf, Ea, Ef, Fa, Fc} can be used to build the five target cubes $Cb, Cd, Ce, Db, De$. The targets $Cb, Cd, Ce$ can be built in four distinct ways, and $Db, De$ can be built in two distinct ways.}
  \label{build5}
\end{figure} 

In the example in Figure~\ref{build5}, we select columns $P = \{ a, c, f\}$. Rows A and F are selected from  $P$, and rows B and E are not in $P$. The five targets, {\it Cb, Cd, Ce, Db,} and {\it De}, are located in the complement of the row and column union.  We remove cubes Bc and Ec from the 10 cubes in the intersection of the selected rows and columns because they are mirror cubes of Cb and Ce; thus, 8 cubes remain for our collection.  Our code verifies that this collection can build  Cb, Cd, and Ce, each with 4 solutions, and build Db and De, each with 2 solutions.

\section{Minimum Universal sets and Proof of Theorem 3}\label{sec:minuniversal}

As we've seen, collections of eight cubes can build no more than five targets, but as mentioned in Section~\ref{sec:Results}, Haraguchi's Minimal Universal Set consists of 12 cubes from which all 30 cubes can be built. How does this transition happen? Can cubes be removed from Haraguchi's set of 12 to give a collection of eight that can build five targets? We will see that the answer to this question is no. What is the maximum number of target cubes that can be built for sets of size 9, 10, or 11?  When, if ever, is it possible to build 29 target cubes? These are both open questions.

We give partial answers by calculating the distribution of the number of targets that can be built for sets of MacMahon cubes of size $k=9, 10, 11,$ and 12, using 20,000 randomly selected sets of each size.  

\begin{figure}[h]
   \centering    \includegraphics[scale=0.22]{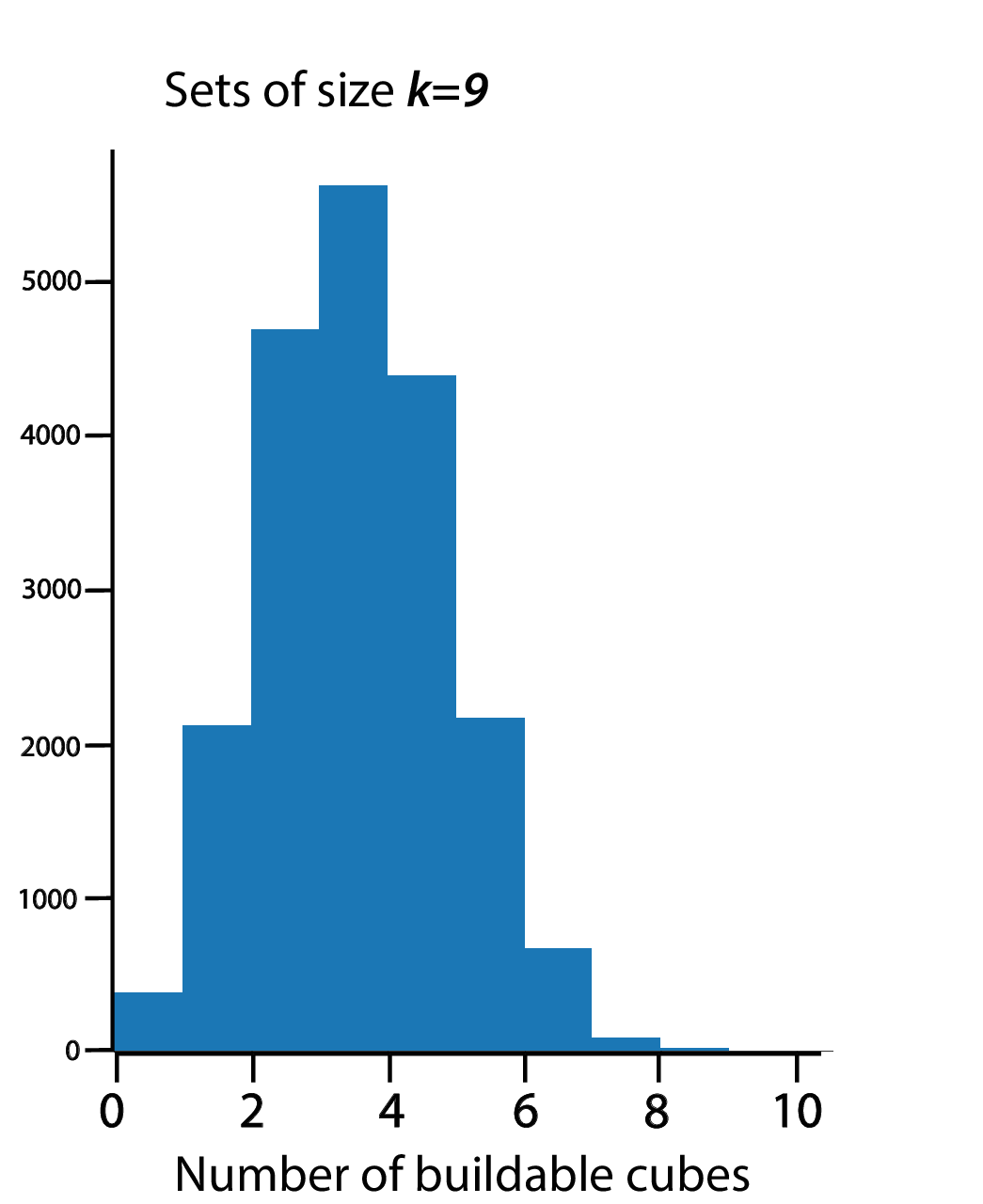} 
   \hspace{0.55in}
   \includegraphics[scale=0.21]{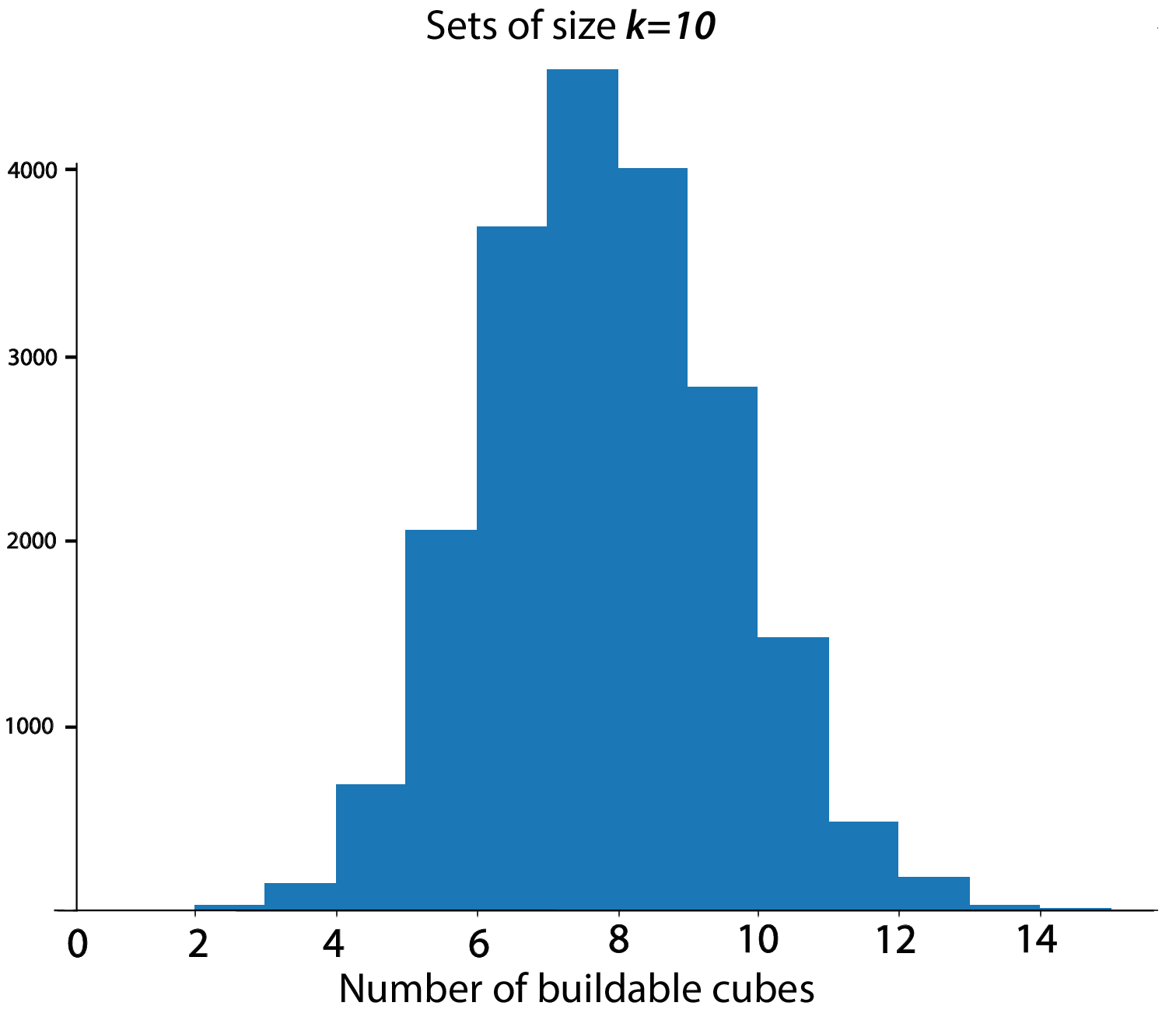} 
   
   \vspace{0.4in}
   
   \includegraphics[scale=0.2]{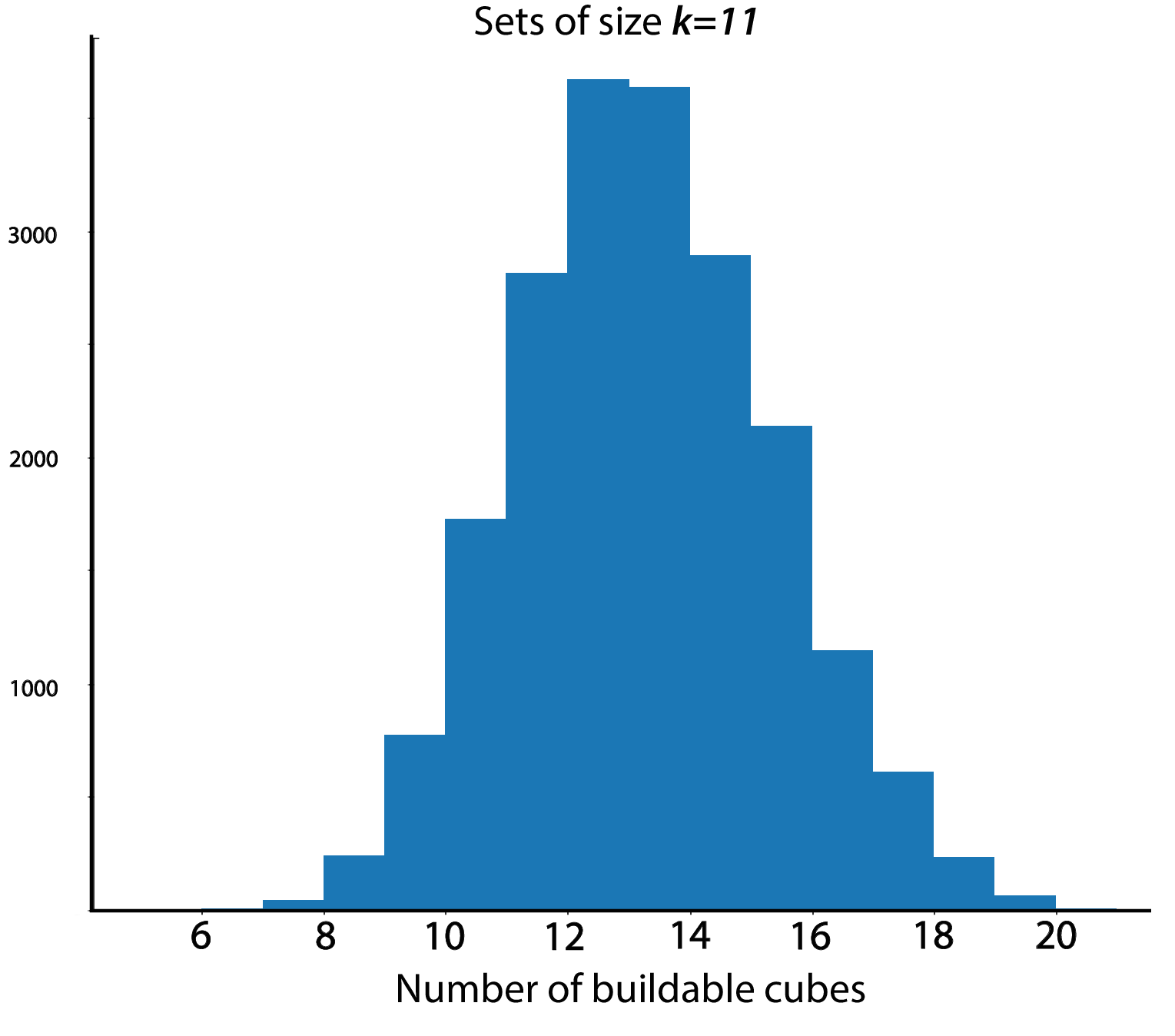}
      \hspace{0.3in}
    \includegraphics[scale=0.2]{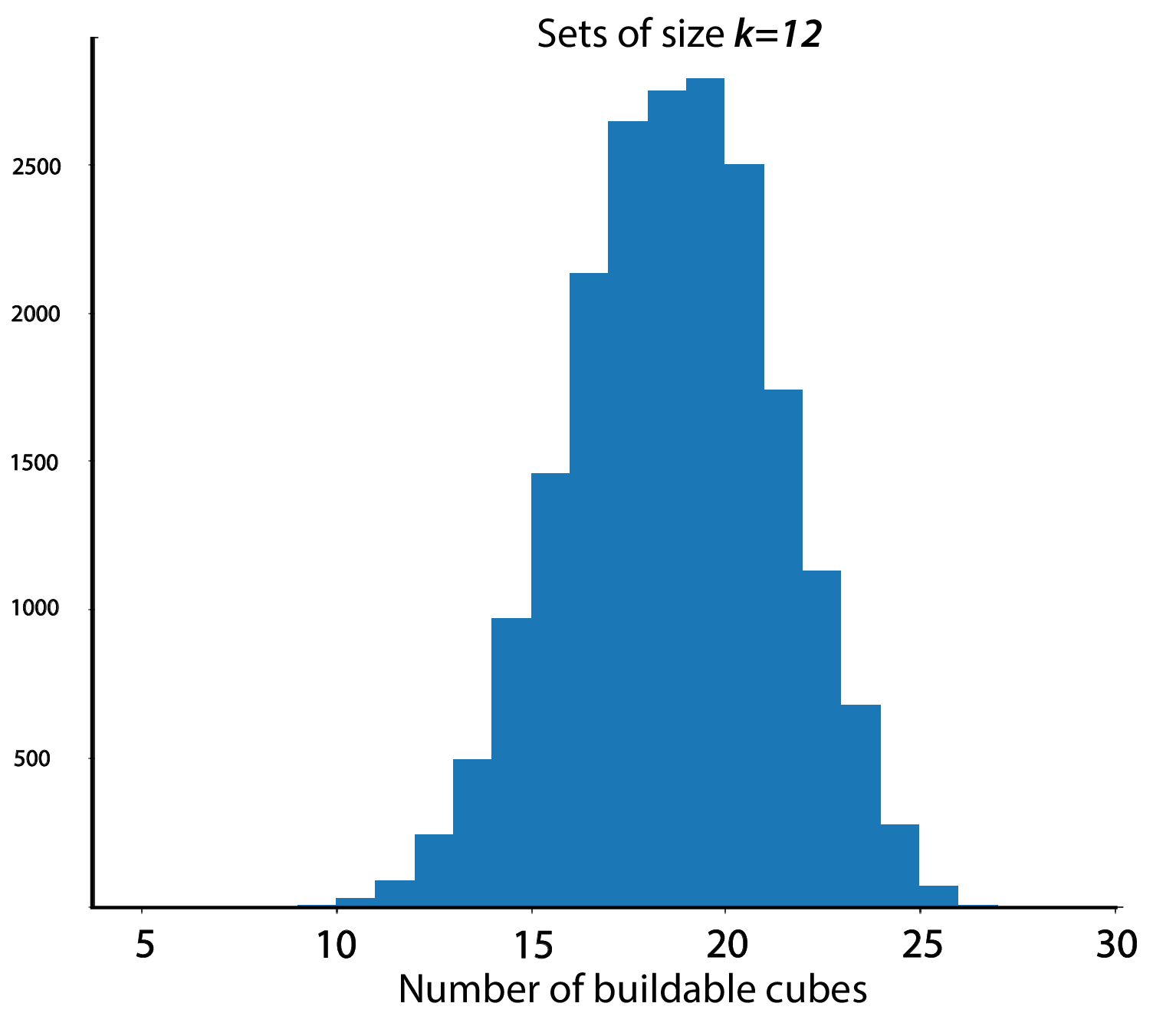}
   \caption{Distributions of the Number of Target Cubes that can be built from 20000 random samples of sets of size $k=9, 10, 11,$ and $12$.}
  \label{20Kdistrib}
\end{figure}

Our calculations show that for 20,000 random sets of nine cubes, the average number of targets that can be built is 3, with a standard deviation of 1.34.  In this calculation for sets of nine cubes a minimum of 0 and a maximum of 8 targets can be built.  For 20,000 random sets of ten cubes, the average number of targets that can be built is 7.3, with a standard deviation of 1.7. For sets of ten cubes a minimum of 2 and a maximum of 14 targets is attained.  Berkove and collaborators prove that for any collection of ten cubes, there is at least one target that can be built~\cite{Berkove2017AutomorphismsOS} and the results here agree with this result.
For 20,000 random sets of eleven cubes, the average number of targets that can be built is  12.8, with a standard deviation of  2.1. For sets of  eleven cubes a minimum of 6 targets and a maximum of 21 targets built is attained.  For 20,000 random sets of twelve cubes, the average number of targets that can be built is 18.2, with a standard deviation of 2.7. For sets of twelve cubes a minimum of 9 and a maximum of 28 targets is attained; but here we know that a maximum of 30 exists. The universal sets of 12 cubes are so rare, they were never encountered in several samples of 20,000 random sets of 12 cubes. Also, we have not yet found any sets of 12 cubes that can build exactly 29 cubes.

Next, we consider Haraguchi's minimum universal set.
\begin{figure}[h]
   \centering    \includegraphics[scale=0.2]{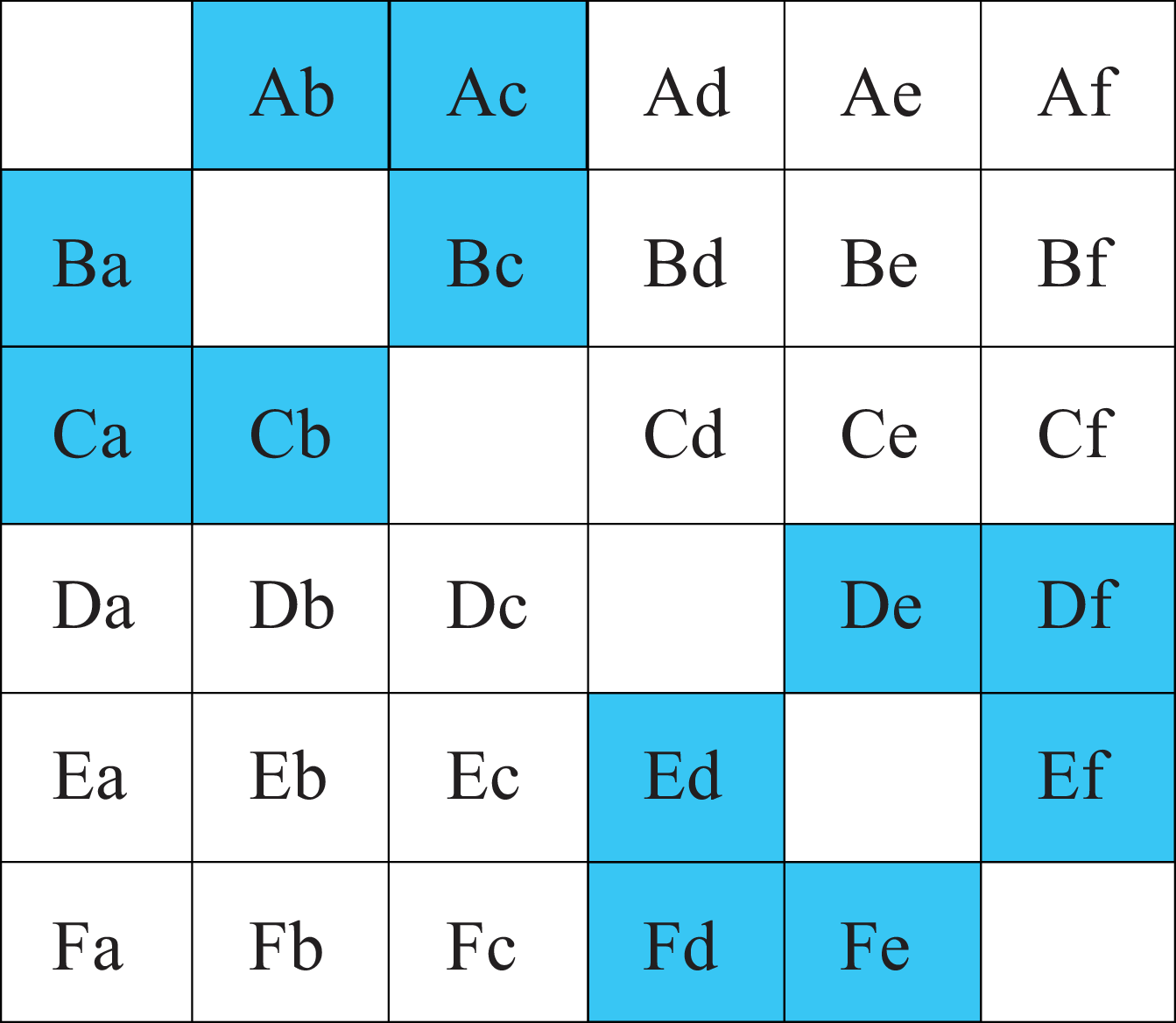}
   \caption{Haraguchi's Minimum Universal set of twelve cubes from which which every target cube can be built.}
  \label{HaraguchiSet}
\end{figure}

From Haraguchi's collection of 12 cubes, we calculate the solution numbers for all 30 target cubes and notice there are two scenarios, depending on whether or not the target is in Haraguchi's set. For a target cube that is in Haraguchi's set, exactly three other cubes from his set have no corner numbers in common with the target cube. These three cubes include the target's mirror cube, the cube in the same row as the target, and the cube in the same column as the target.  Excluding these three cubes the remaining nine cubes give nine possible collections to build the target.  Our calculations show that two of the nine collections solve the puzzle with a solution number of 8, and the remaining seven collections also solve the puzzle but with a solution number 2. 
For each of the 18 cubes that are not in Haraguchi's set, there are exactly four Haraguchi cubes that lie in the same row or column as the target and these cubes have no corner numbers in common with the target cube.  Since these four Haraguchi cubes are useless in this target puzzle, there is only one remaining collection to build the target, and that collection can always solve the puzzle with a solution number of 4.

Another question to consider is that Haraguchi's set of 12 cubes can build all 30 cubes, while an arbitrary collection of 8 cubes can build no more than 5 cubes.  For Haraguchi's cubes, we investigate how the distribution of the number of targets that can be built changes if we decrease from all 12 cubes to sets of 11, 10, 9, or 8 cubes.  Are there any collections of 8 cubes in Haraguchi's set that build the maximum of 5 cubes?  Our calculations show that the answer is no. All collections from Haraguchi's set can build no more than 3 cubes.

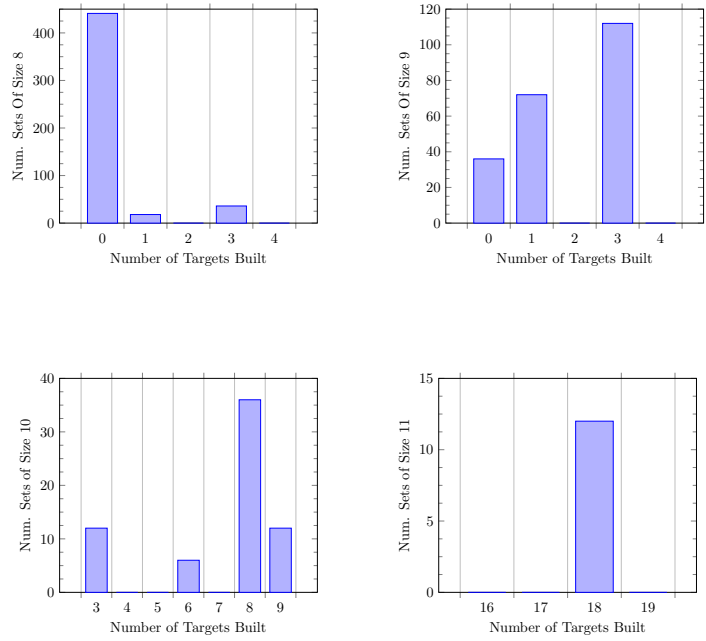
\begin{figure}
    \centering
    \begin{tikzpicture}[scale=0.5]
\begin{axis}[ybar interval=0.7, ymax=450,ymin=0, minor y tick num = 3,  ylabel=Num. Sets Of Size 8, xlabel= Number of Targets Built]
\addplot coordinates { (0, 441) (1, 18) (2, 0) (3, 36) (4, 0) (5,0)};
\end{axis}
\end{tikzpicture}
\hspace{0.3in}
    \begin{tikzpicture}[scale=0.5]
\begin{axis}[ybar interval=0.7, ymax=120, ymin=0, minor y tick num = 3, ylabel= Num. Sets Of Size 9, xlabel= Number of Targets Built]
\addplot coordinates { (0, 36) (1, 72) (2, 0) (3, 112) (4, 0) (5,0) };
\end{axis}
\end{tikzpicture}

\vspace{0.5in}

    \begin{tikzpicture}[scale=0.5]
\begin{axis}[ybar interval=0.7, ymax=40,ymin=0, minor y tick num = 3,  ylabel=Num. Sets of Size 10, xlabel= Number of Targets Built]
\addplot coordinates {(3, 12) (4,0) (5,0) (6, 6) (7,0) (8, 36) (9, 12) (10, 0)};
\end{axis}
\end{tikzpicture}
\hspace{0.3in}
    \centering
    \begin{tikzpicture}[scale=0.5]
\begin{axis}[ybar interval=0.7, ymax=15,ymin=0, minor y tick num = 3, ylabel=Num. Sets of Size 11, xlabel= Number of Targets Built]
\addplot coordinates {(16,0) (17,0) (18, 12) (19, 0) (20, 0)};
\end{axis}
\end{tikzpicture}
    \caption{The distribution of the number of target cubes built from subsets of size $k \in \{8, 9, 10, 11\}$ from Haraguchi's set.}
   \label{DistBuildHaraguchi}
\end{figure}

In Figure~\ref{DistBuildHaraguchi},  we see that subsets of eight Haraguchi cubes build no more than 3 targets, and most build no targets at all.  When considering subsets of nine Haraguchi cubes,  many build 3 targets, but some sets build none. For subsets of ten Haraguchi cubes, the number of targets built ranges from 3 to 9 targets, but they never build 4, 5, or 7 targets. All subsets of eleven Haraguchi cubes build 18 targets. When we consider the set of all 12 Haraguchi cubes,  the number of targets that can be built increases from 18 to all 30 cubes. 

\begin{figure}[h]
   \centering  \includegraphics[scale=0.15]{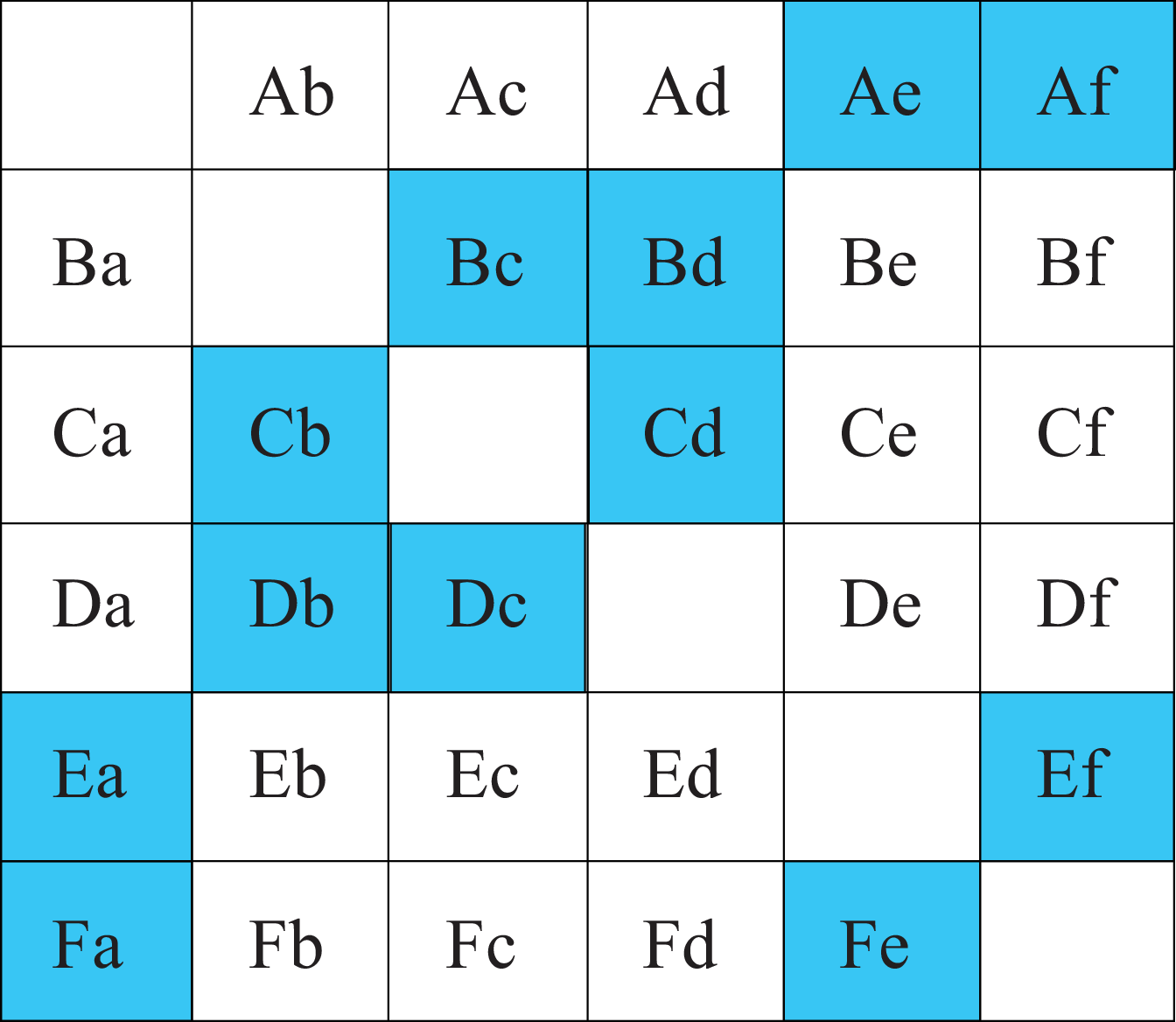} 
   \hspace{0.3in} \includegraphics[scale=0.15]{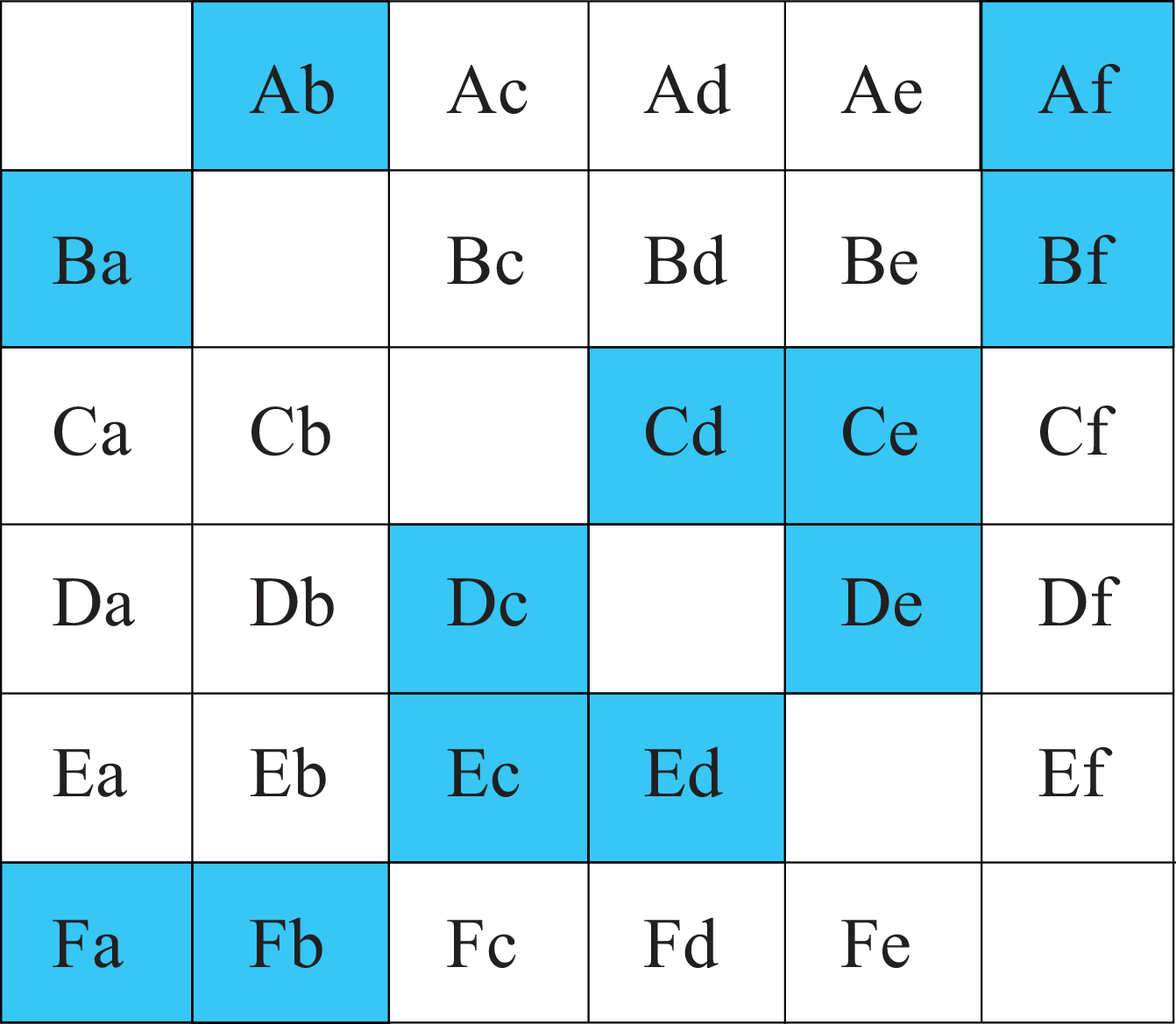} 
   \caption{Two additional minimum universal sets from which the target puzzle can be solved for all target cubes. These two collections are versions of Haraguchi's set but shifted along the tableau's diagonal, assuming torus identifications on the perimeter of the tableau. }
  \label{TwoMore}
\end{figure}

\begin{figure}[h]
   \centering  \includegraphics[scale=0.15]{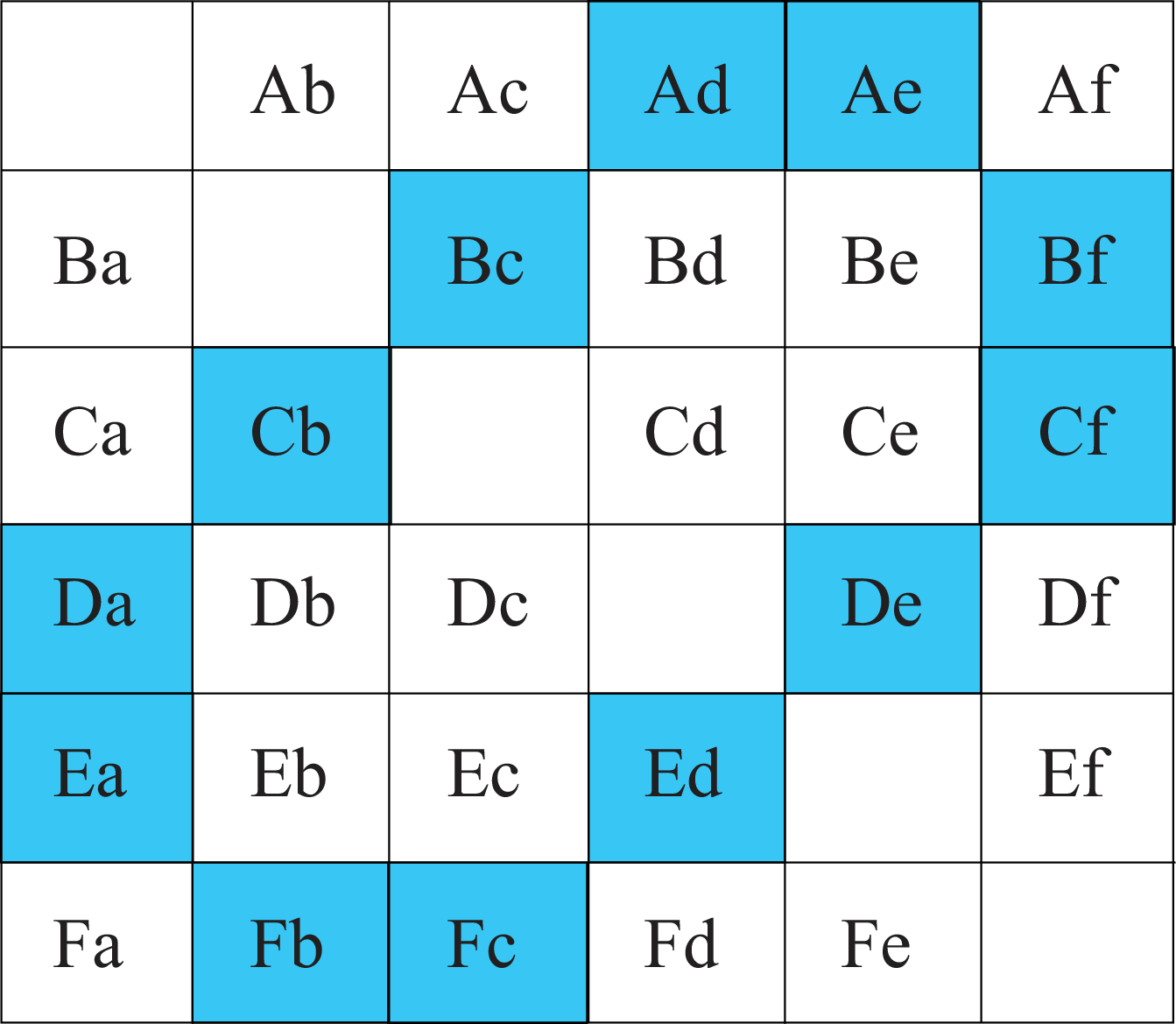} 
   \hspace{0.3in} 
   \includegraphics[scale=0.15]{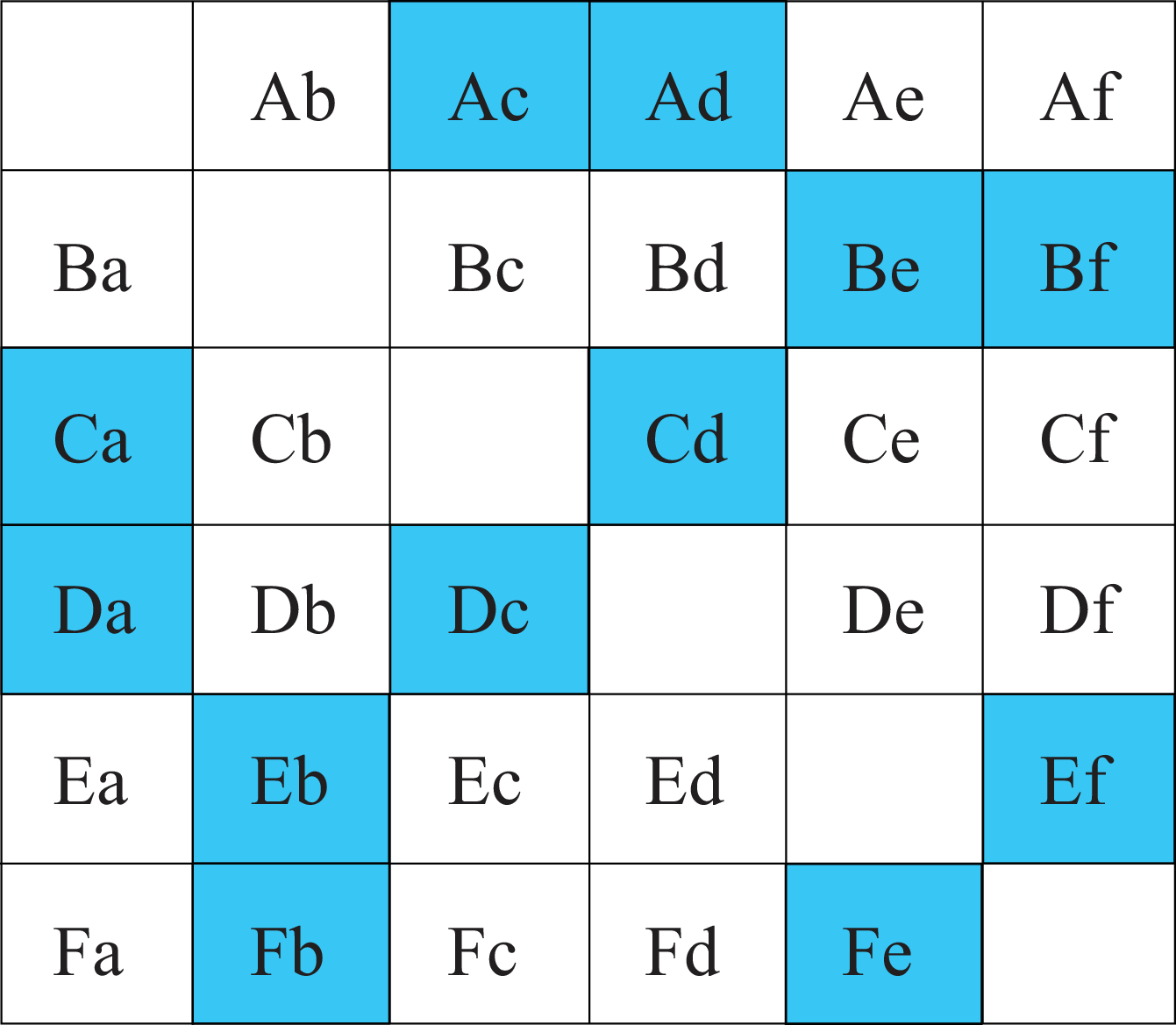} 
    \hspace{0.3in}
   \includegraphics[scale=0.15]{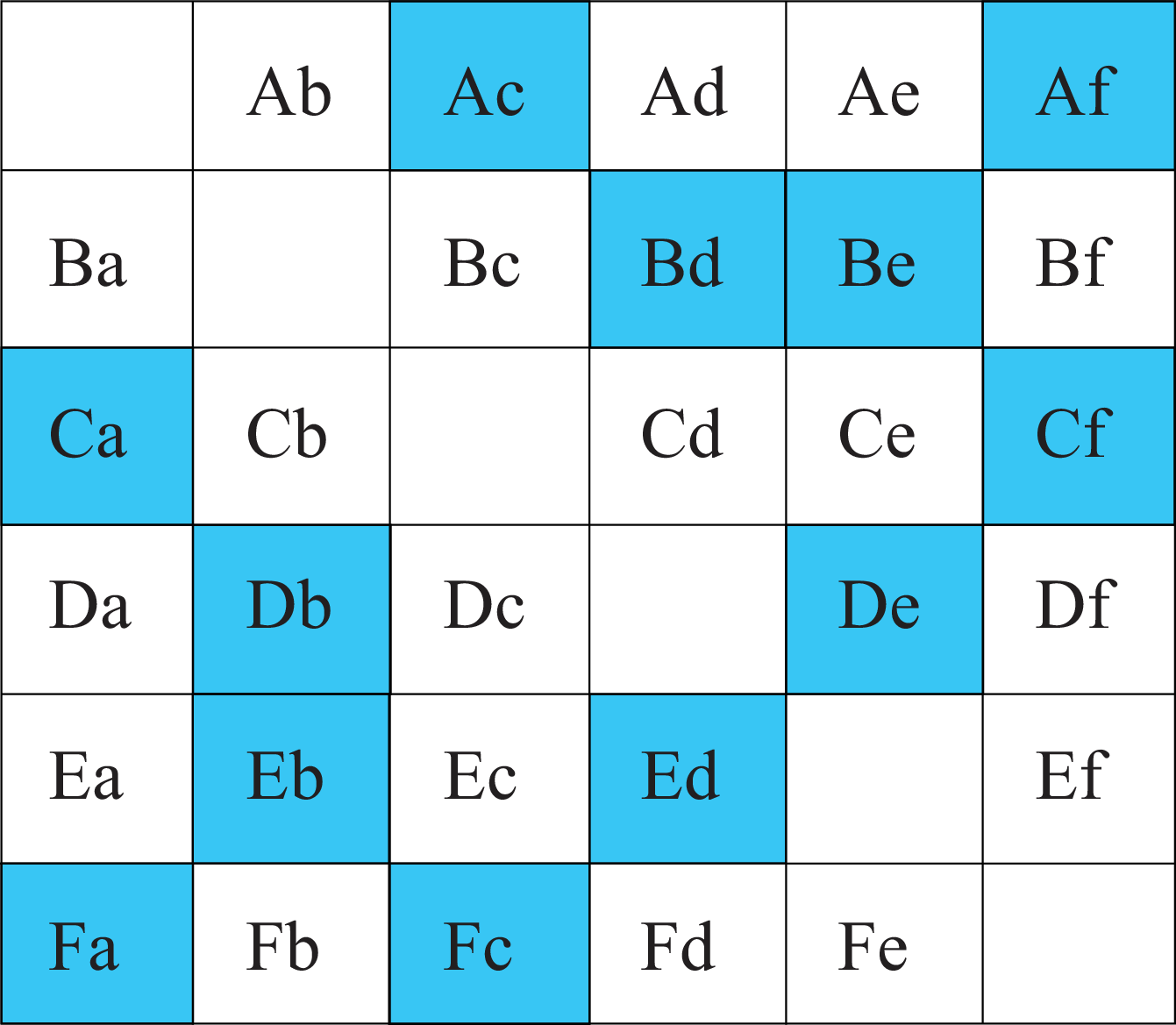}

    \
    
  \includegraphics[scale=0.15]{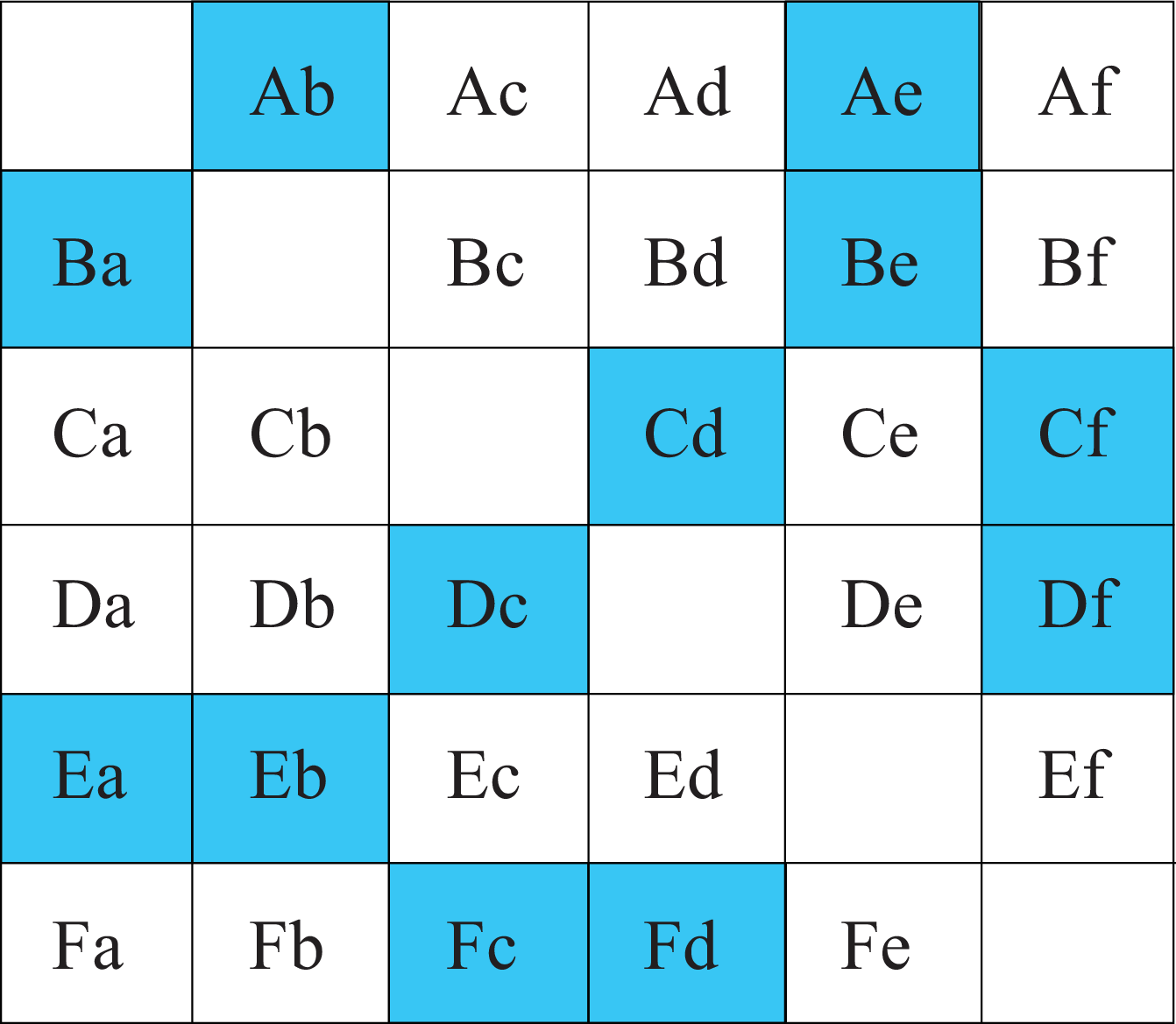} 
   \hspace{0.3in}
  \includegraphics[scale=0.15]{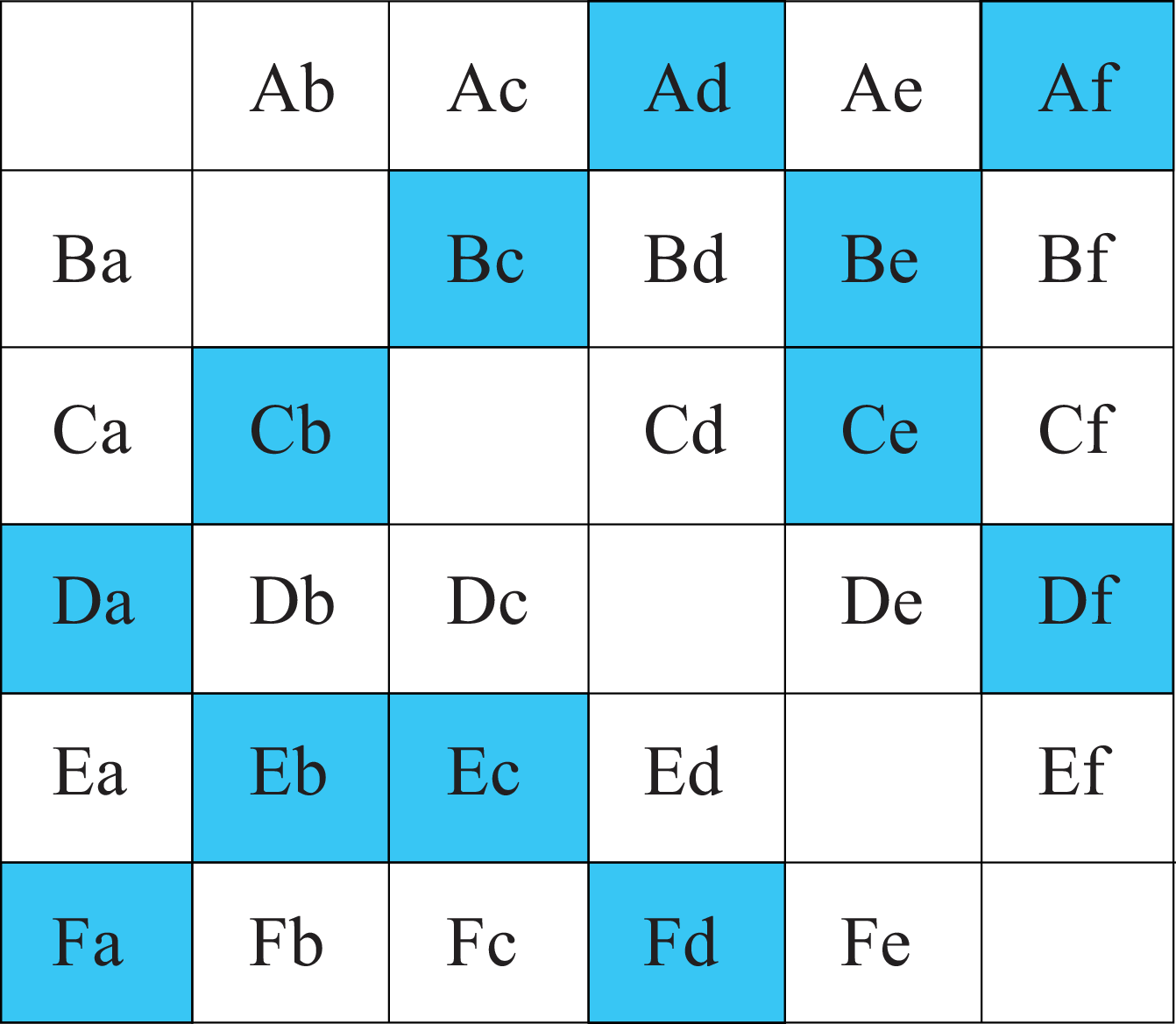} 
  \hspace{0.3in}
    \includegraphics[scale=0.15]{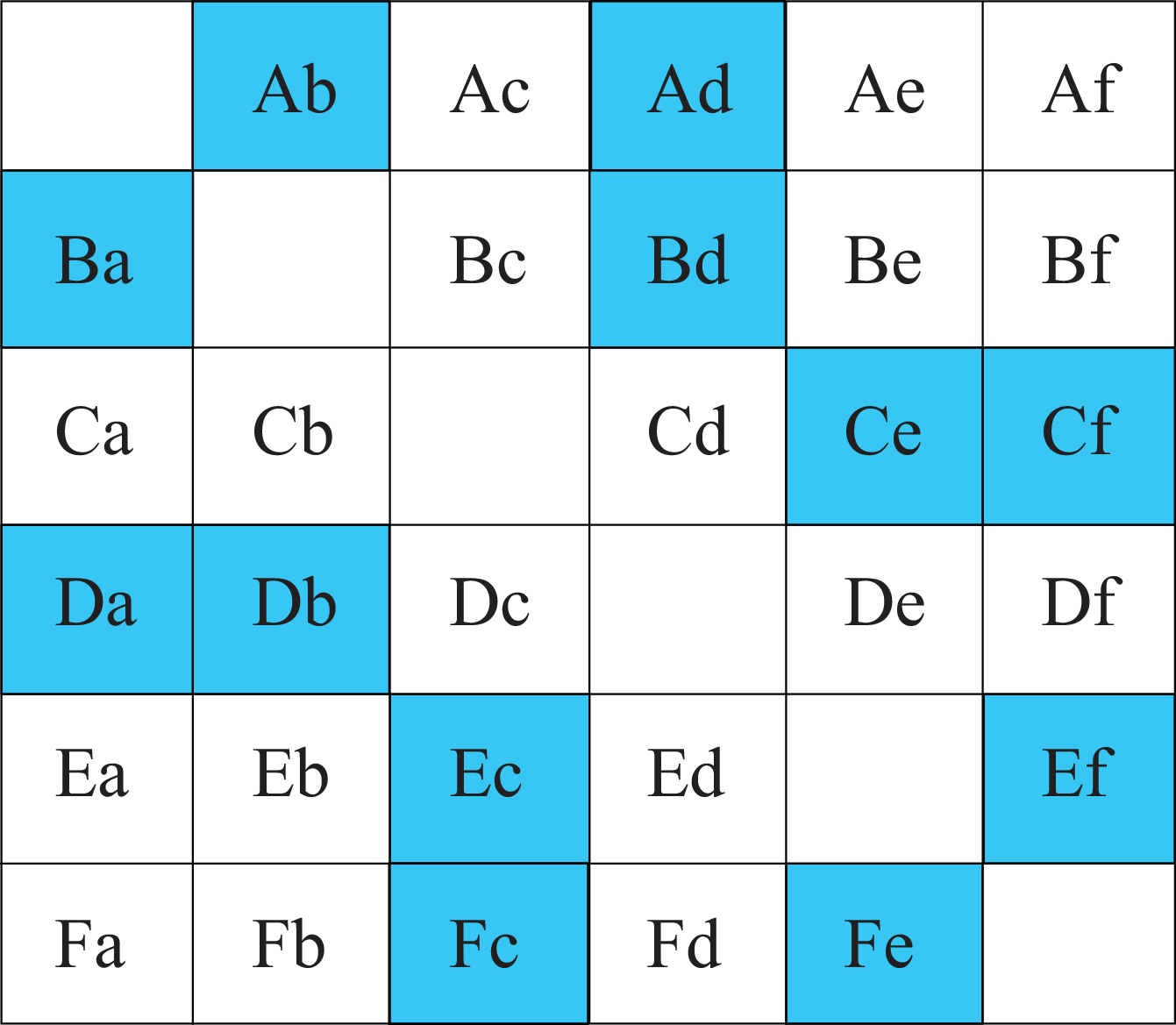} 
   \caption{Six more minimum universal sets from which the target puzzle can be solved for all 30 cubes in $C_6$.}
  \label{SixMore}
\end{figure}

We found two shifted versions of Haraguchi's set of 12 cubes that can also be used to solve the target puzzle for every MacMahon cube.  These two sets are shown in Figure~\ref{TwoMore}. Notice that they are shifted versions of Haraguchi's set along the diagonal, assuming torus identifications along the edges of the tableau.  We found another six minimum universal sets of 12 cubes; all such sets are symmetric across the diagonal and consist of two cubes in each column and each row. These six sets, shown in Figure~\ref{SixMore}, are also shifted versions of each other down the diagonal, assuming torus identifications along the edges of the tableau. 

There are many collections of 12 cubes that are symmetric along the diagonal and contain two cubes from each row and column, but very few are minimum universal sets.  The following theorem describes our conjecture for the special sets of 12 cubes that can build all 30 MacMahon cubes.

\begin{conjecture}\label{MinUniversalT}
    Every minimum universal set, $M$, of 12 cubes satisfies the following rules.
\begin{itemize}
    \item $M$ contains two cubes from Ab, Ac, Ad, Ae, Af. We refer to these cubes as Ax and Ay for some x, y $\in \{$a, b, c, d, e, f $\}$.  
    \item $M$ also contains Xa, Ya, Xy, and Yx. 
    \item The remaining six cubes are the complementary letter pairs to the six specified above. That is, $M$ contains Wz, Wv, Zw, Zv, Vw, Vz, where the variables W, Z, V, X, Y, and A comprise distinct elements in the set $\{$A, B, C, D, E, F$\}$.
\end{itemize}
\end{conjecture}

\begin{figure}[h]
   \centering  \includegraphics[scale=0.15]{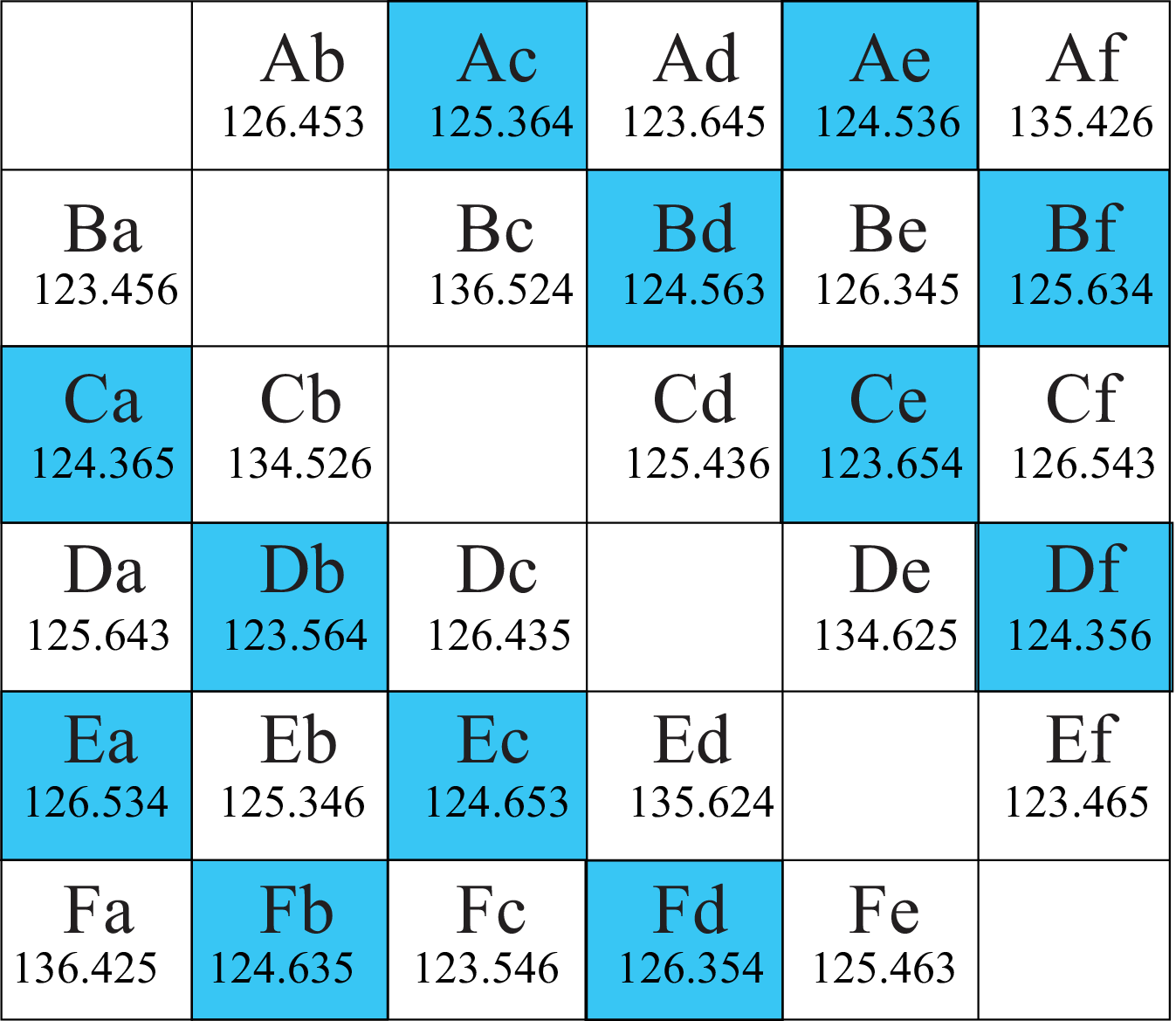} 
   \caption{The tenth minimum universal set for the 30 MacMahon cubes.}
  \label{noshifts}
\end{figure}

Notice that Conjecture~\ref{MinUniversalT} implies that there are exactly ten minimum universal sets because there are 10 ways to select two cubes from Ab, Ac, Ad, Ae, Af. We have seen nine minimum universal sets above, and the tenth is shown in Figure~\ref{noshifts}.  Notice that any shift of this collection along the diagonal, assuming torus identifications along the edges, results in the same collection.

These ten minimum universal sets of 12 cubes are all related by permuting the cube colors.  In other words, they form an orbit of the $S_6$ group action on sets of 12 McMahon cubes. We have observed many 3-cycles within this orbit that send a minimum universal set of cubes to itself. However, these 3-cycles are not the identity on individual cubes.


\section{Algorithm for Computations}\label{sec:alg}


The insights found in Theorems~\ref{WnoTarget} and \ref{WhasTarget} arose from investigations using an algorithm that calculates the solution number for a given target over a family of collections of eight cubes. In this section, we describe this algorithm and provide a brief example calculation.

{\bf Input:} target cube $T \in C_k$ and a dictionary, $\mathcal{C}$, of cubes in $C_k$ with (cube name, \{corner numbers\}) as the (key, \{value\}) pairs. Usually, $\mathcal{C} = C_k$, but we can also restrict to collections from a subset of $C_k$.
\begin{enumerate}
\item {\bf Create a subset of usable cubes for $T$.}
Start by creating a collection of usable cubes for our target cube, denoted $U_T$, by removing cubes with a usable corner count of zero. $$U_T:= \{ Y \in \mathcal{C} \ |  \ ucc(T, Y) \geq 1\}.$$

\item {\bf Generate a list $
\mathcal{L}$ of all collections of 8 distinct cubes from $U_T$.} \newline The following steps are applied to each
collection $C$ in $\mathcal{L}.$

\begin{enumerate}
\item {\bf Create the Corner Count Vector for a given collection.}
Let $C$ be a collection of 8 distinct cubes from $U_T$. We associate the ordered list of 8 cubes with an index from a list of the first 8 primes 
$$[2, 3, 5, 7, 11, 13, 17, 19].$$
In other words, we will refer to the first cube in the list of 8 as having `index' 2, and the fourth cube in the list of 8 as having `index' 7, and so on.

For each corner $x_i$ of $T$, we create the $i^{th}$ entry of our corner count vector as follows
\begin{itemize}
\item $m_{i}$ is the number of cubes in $C$ that contain corner $x_i$, and 
\item $[p_{i,1}, \dots, p_{i,m_{i}}]$ is a list of prime indices of the cubes that contain corner $x_i$.
\end{itemize}
Then the corner count vector is $$[[x_1, m_{1}, [p_{1,1}, \dots, p_{1,m_{x_1}}]], \dots ,
[x_8, m_{8}, [p_{8,1}, \dots, p_{8,m_{x_8}}]].$$

If $m_i = 0$ for some corner $x_i$ of $T$, then the number of solutions for the collection $C$ is zero. {\bf Record a solution of 0 and stop algorithm for this collection $C$.}

\item {\bf Determine the number of solutions for $T$ given $C$.}
Create the Cartesian product of the eight nonempty sets of primes, $$CP = \Pi_{i=1}^8 \{p_{i,1}, \dots, p_{i,m_{i}}\}.$$ 
Each tuple $(q_1, q_2, \dots , q_8)\in CP$, is viewed as an {\it attempted} solution of the target puzzle for $T$ with the cube of index $q_i$ placed in corner $x_i$ for $i=1, 2, \dots, 8$. 
We check the tuple to see if it contributes to the solution count for $T$ using the following congruence.  If 
$$ q_1 \cdot q_2 \dots \cdot q_8 \equiv 0 \pmod{ 2 \cdot 3 \cdot 5  \dots \cdot 19} $$
then every cube in the collection is used and every corner of $T$  is satisfied by some cube in the collection. 
{\bf Count the number of tuples that result in a solution for $T$ given $C$}.
\end{enumerate}

\item {\bf Store, sort and display all the non-zero solution number(s) for $T$ for all the collections in $\mathcal{L}$.}  
\end{enumerate}

{\bf Output:} The distribution of solution numbers for the target puzzle for $T$ using collections of eight cubes from the dictionary $\mathcal{C}$.


\normalsize
\begin{table}[h]
\begin{tabular}{|c|llllllll|}
\hline
\multicolumn{1}{|c|}{Name} &
  \multicolumn{8}{c|}{Corners} \\ \hline
Ab &
  \multicolumn{1}{l|}{143} &
  \multicolumn{1}{l|}{345} &
  \multicolumn{1}{l|}{235} &
  \multicolumn{1}{l|}{132} &
  \multicolumn{1}{l|}{126} &
  \multicolumn{1}{l|}{256} &
  \multicolumn{1}{l|}{465} &
  164 \\ \hline
Ac &
  \multicolumn{1}{l|}{153} &
  \multicolumn{1}{l|}{134} &
  \multicolumn{1}{l|}{142} &
  \multicolumn{1}{l|}{125} &
  \multicolumn{1}{l|}{265} &
  \multicolumn{1}{l|}{246} &
  \multicolumn{1}{l|}{364} &
  356 \\ \hline
Ad &
  \multicolumn{1}{l|}{243} &
  \multicolumn{1}{l|}{123} &
  \multicolumn{1}{l|}{152} &
  \multicolumn{1}{l|}{254} &
  \multicolumn{1}{l|}{456} &
  \multicolumn{1}{l|}{165} &
  \multicolumn{1}{l|}{136} &
  346 \\ \hline
Ae &
  \multicolumn{1}{l|}{354} &
  \multicolumn{1}{l|}{145} &
  \multicolumn{1}{l|}{124} &
  \multicolumn{1}{l|}{234} &
  \multicolumn{1}{l|}{263} &
  \multicolumn{1}{l|}{162} &
  \multicolumn{1}{l|}{156} &
  365 \\ \hline
Af &
  \multicolumn{1}{l|}{245} &
  \multicolumn{1}{l|}{154} &
  \multicolumn{1}{l|}{135} &
  \multicolumn{1}{l|}{253} &
  \multicolumn{1}{l|}{236} &
  \multicolumn{1}{l|}{163} &
  \multicolumn{1}{l|}{146} &
  264 \\ \hline
  
Ba &
  \multicolumn{1}{l|}{123} &
  \multicolumn{1}{l|}{253} &
  \multicolumn{1}{l|}{354} &
  \multicolumn{1}{l|}{134} &
  \multicolumn{1}{l|}{146} &
  \multicolumn{1}{l|}{456} &
  \multicolumn{1}{l|}{265} &
  162 \\ \hline
Bc &
  \multicolumn{1}{l|}{143} &
  \multicolumn{1}{l|}{154} &
  \multicolumn{1}{l|}{245} &
  \multicolumn{1}{l|}{234} &
  \multicolumn{1}{l|}{263} &
  \multicolumn{1}{l|}{256} &
  \multicolumn{1}{l|}{165} &
  136 \\ \hline
Bd &
  \multicolumn{1}{l|}{153} &
  \multicolumn{1}{l|}{132} &
  \multicolumn{1}{l|}{124} &
  \multicolumn{1}{l|}{145} &
  \multicolumn{1}{l|}{465} &
  \multicolumn{1}{l|}{264} &
  \multicolumn{1}{l|}{236} &
  356 \\ \hline
Be &
  \multicolumn{1}{l|}{152} &
  \multicolumn{1}{l|}{135} &
  \multicolumn{1}{l|}{345} &
  \multicolumn{1}{l|}{254} &
  \multicolumn{1}{l|}{246} &
  \multicolumn{1}{l|}{364} &
  \multicolumn{1}{l|}{163} &
  126 \\ \hline
Bf &
  \multicolumn{1}{l|}{243} &
  \multicolumn{1}{l|}{235} &
  \multicolumn{1}{l|}{125} &
  \multicolumn{1}{l|}{142} &
  \multicolumn{1}{l|}{164} &
  \multicolumn{1}{l|}{156} &
  \multicolumn{1}{l|}{365} &
  346 \\ \hline
Ca &
  \multicolumn{1}{l|}{152} &
  \multicolumn{1}{l|}{124} &
  \multicolumn{1}{l|}{143} &
  \multicolumn{1}{l|}{135} &
  \multicolumn{1}{l|}{365} &
  \multicolumn{1}{l|}{346} &
  \multicolumn{1}{l|}{264} &
  256 \\ \hline

Cb &
  \multicolumn{1}{l|}{243} &
  \multicolumn{1}{l|}{254} &
  \multicolumn{1}{l|}{145} &
  \multicolumn{1}{l|}{134} &
  \multicolumn{1}{l|}{163} &
  \multicolumn{1}{l|}{156} &
  \multicolumn{1}{l|}{265} &
  236 \\ \hline
Cd &
  \multicolumn{1}{l|}{235} &
  \multicolumn{1}{l|}{345} &
  \multicolumn{1}{l|}{154} &
  \multicolumn{1}{l|}{125} &
  \multicolumn{1}{l|}{162} &
  \multicolumn{1}{l|}{146} &
  \multicolumn{1}{l|}{364} &
  263 \\ \hline
Ce &
  \multicolumn{1}{l|}{245} &
  \multicolumn{1}{l|}{253} &
  \multicolumn{1}{l|}{123} &
  \multicolumn{1}{l|}{142} &
  \multicolumn{1}{l|}{164} &
  \multicolumn{1}{l|}{136} &
  \multicolumn{1}{l|}{356} &
  465 \\ \hline
Cf &
  \multicolumn{1}{l|}{153} &
  \multicolumn{1}{l|}{354} &
  \multicolumn{1}{l|}{234} &
  \multicolumn{1}{l|}{132} &
  \multicolumn{1}{l|}{126} &
  \multicolumn{1}{l|}{246} &
  \multicolumn{1}{l|}{456} &
  165 \\ \hline
Da &
  \multicolumn{1}{l|}{245} &
  \multicolumn{1}{l|}{125} &
  \multicolumn{1}{l|}{132} &
  \multicolumn{1}{l|}{234} &
  \multicolumn{1}{l|}{364} &
  \multicolumn{1}{l|}{163} &
  \multicolumn{1}{l|}{156} &
  465 \\ \hline
Db &
  \multicolumn{1}{l|}{154} &
  \multicolumn{1}{l|}{142} &
  \multicolumn{1}{l|}{123} &
  \multicolumn{1}{l|}{135} &
  \multicolumn{1}{l|}{365} &
  \multicolumn{1}{l|}{263} &
  \multicolumn{1}{l|}{246} &
  456 \\ \hline
Dc &
  \multicolumn{1}{l|}{152} &
  \multicolumn{1}{l|}{145} &
  \multicolumn{1}{l|}{354} &
  \multicolumn{1}{l|}{253} &
  \multicolumn{1}{l|}{236} &
  \multicolumn{1}{l|}{346} &
  \multicolumn{1}{l|}{164} &
  126 \\ \hline
De &
  \multicolumn{1}{l|}{153} &
  \multicolumn{1}{l|}{235} &
  \multicolumn{1}{l|}{243} &
  \multicolumn{1}{l|}{134} &
  \multicolumn{1}{l|}{146} &
  \multicolumn{1}{l|}{264} &
  \multicolumn{1}{l|}{256} &
  165 \\ \hline
Df &
  \multicolumn{1}{l|}{143} &
  \multicolumn{1}{l|}{124} &
  \multicolumn{1}{l|}{254} &
  \multicolumn{1}{l|}{345} &
  \multicolumn{1}{l|}{356} &
  \multicolumn{1}{l|}{265} &
  \multicolumn{1}{l|}{162} &
  136 \\ \hline
Ea &
  \multicolumn{1}{l|}{243} &
  \multicolumn{1}{l|}{142} &
  \multicolumn{1}{l|}{154} &
  \multicolumn{1}{l|}{345} &
  \multicolumn{1}{l|}{356} &
  \multicolumn{1}{l|}{165} &
  \multicolumn{1}{l|}{126} &
  236 \\ \hline
Eb &
  \multicolumn{1}{l|}{245} &
  \multicolumn{1}{l|}{354} &
  \multicolumn{1}{l|}{153} &
  \multicolumn{1}{l|}{125} &
  \multicolumn{1}{l|}{162} &
  \multicolumn{1}{l|}{136} &
  \multicolumn{1}{l|}{346} &
  264 \\ \hline 
Ec &
  \multicolumn{1}{l|}{124} &
  \multicolumn{1}{l|}{132} &
  \multicolumn{1}{l|}{235} &
  \multicolumn{1}{l|}{254} &
  \multicolumn{1}{l|}{456} &
  \multicolumn{1}{l|}{365} &
  \multicolumn{1}{l|}{163} &
  146 \\ \hline
Ed &
  \multicolumn{1}{l|}{143} &
  \multicolumn{1}{l|}{234} &
  \multicolumn{1}{l|}{253} &
  \multicolumn{1}{l|}{135} &
  \multicolumn{1}{l|}{156} &
  \multicolumn{1}{l|}{265} &
  \multicolumn{1}{l|}{246} &
  164 \\ \hline
Ef &
  \multicolumn{1}{l|}{152} &
  \multicolumn{1}{l|}{123} &
  \multicolumn{1}{l|}{134} &
  \multicolumn{1}{l|}{145} &
  \multicolumn{1}{l|}{465} &
  \multicolumn{1}{l|}{364} &
  \multicolumn{1}{l|}{263} &
  256 \\ \hline
Fa &
  \multicolumn{1}{l|}{235} &
  \multicolumn{1}{l|}{153} &
  \multicolumn{1}{l|}{145} &
  \multicolumn{1}{l|}{254} &
  \multicolumn{1}{l|}{246} &
  \multicolumn{1}{l|}{164} &
  \multicolumn{1}{l|}{136} &
  263 \\ \hline
Fb &
  \multicolumn{1}{l|}{124} &
  \multicolumn{1}{l|}{152} &
  \multicolumn{1}{l|}{253} &
  \multicolumn{1}{l|}{234} &
  \multicolumn{1}{l|}{364} &
  \multicolumn{1}{l|}{356} &
  \multicolumn{1}{l|}{165} &
  146 \\ \hline
Fc &
  \multicolumn{1}{l|}{123} &
  \multicolumn{1}{l|}{243} &
  \multicolumn{1}{l|}{345} &
  \multicolumn{1}{l|}{135} &
  \multicolumn{1}{l|}{156} &
  \multicolumn{1}{l|}{465} &
  \multicolumn{1}{l|}{264} &
  162 \\ \hline
Fd &
  \multicolumn{1}{l|}{354} &
  \multicolumn{1}{l|}{245} &
  \multicolumn{1}{l|}{142} &
  \multicolumn{1}{l|}{134} &
  \multicolumn{1}{l|}{163} &
  \multicolumn{1}{l|}{126} &
  \multicolumn{1}{l|}{256} &
  365 \\ \hline
Fe &
  \multicolumn{1}{l|}{154} &
  \multicolumn{1}{l|}{143} &
  \multicolumn{1}{l|}{132} &
  \multicolumn{1}{l|}{125} &
  \multicolumn{1}{l|}{265} &
  \multicolumn{1}{l|}{236} &
  \multicolumn{1}{l|}{346} &
  456 \\ \hline
\end{tabular}\caption{Corner numbers for all 30 6-colored cubes.}
    \label{C36corners}
\end{table}

\bibliographystyle{plain}

\bibliography{mainbib}

\end{document}